\documentclass[11pt,a4paper]{article}
\usepackage[T1]{fontenc}
\usepackage{extarrows}
\usepackage{tikz-cd}
\usepackage{lmodern}
\usepackage[utf8]{inputenc}
\usepackage[english]{babel}
\usepackage[margin=28mm]{geometry}
\usepackage{amsmath,amssymb,amsthm,mathtools}
\usepackage{booktabs,tabularx,array}
\usepackage[hidelinks]{hyperref}
\usepackage[expansion=false]{microtype}
\newcommand{\GL}{\mathrm{GL}}

\newcommand{\CE}{d_{\g}}
\newcommand{\nil}{\mathfrak n}

\DeclareMathOperator{\rank}{rank}

\newcommand{\Span}{\mathrm{Span}}

\newcommand{\G}{{\mathfrak{g}}}

\newcommand{\kalg}{\mathfrak k}
\newcommand{\h}{{\mathfrak{h}}}
\newcommand{\s}{{\mathfrak{s}}}

\newcommand{\ad}{\mathrm{ad}}

\newcommand{\esp}{\quad\mbox{and}\quad}
\newcommand{\R}{\mathbb{R}}

\newcommand{\om}{\omega}

\newcommand{\too}{\longrightarrow}

\newcommand{\wi}{\widetilde}

\newcommand{\g}{\mathfrak{g}}
\newcommand{\rf}{\mathfrak{r}}

\newcommand{\Ad}{{\mathrm{Ad}}}

\newcommand{\Om}{\Omega}

\newcommand{\al}{\alpha}
\newcommand{\be}{\beta}
\newcommand{\ga}{\gamma}

\newcommand{\la}{\lambda}

\newtheorem{theorem}{Theorem}[section]
\newtheorem{proposition}{Proposition}[section]

\newtheorem{lemma}[theorem]{Lemma}
\newtheorem{corollary}[theorem]{Corollary}
\theoremstyle{definition}

\theoremstyle{remark}
\newtheorem{remark}[theorem]{Remark}

\numberwithin{equation}{section}
\allowdisplaybreaks
\title{On compact homogeneous Locally conformal symplectic manifolds }
\author{M.~Boucetta}
\date{email: m.boucetta@uca.ac.ma}

\begin{document}
	\maketitle
	
	\begin{abstract} A locally conformally symplectic (LCS) structure on a manifold $M$ is a pair $(\eta,\omega)$ consisting of a closed $1$-form $\eta$ called Lee form and a nondegenerate $2$-form $\omega$ satisfying
		\[
		d\omega=\eta\wedge\omega.
		\]
		
		We classify compact connected homogeneous LCS manifolds with nonzero Lee form. To this end, we construct three explicit families of pairs $(G,K)$, where $G$ is a connected simply connected Lie group and $K$ is a closed subgroup such that $G/K$ is compact, the $G$-action is almost effective, and $G/K$ admits a $G$-invariant LCS structure. We show that every compact connected homogeneous LCS manifold with nonzero Lee form arises from one of these families.
		
		Moreover, each such manifold is LCS-isomorphic to a product $(S/H_d)\times S^1$, where $S$ is a compact connected simply connected semisimple Lie group and $S/H_d$ carries an $S$-invariant contact form $\alpha$. Under this identification, the LCS structure is given by
		\[
		\eta=\pi_2^*\xi,\qquad
		\omega=d(\pi_1^*\alpha)+\pi_1^*\alpha\wedge\eta,
		\]
		where $\pi_1$ and $\pi_2$ are the projections onto the two factors and $\xi$ is a nowhere-vanishing $1$-form on $S^1$. This description also yields topological consequences for compact homogeneous LCS manifolds. In particular,
		none of these compact
		homogeneous LCS manifolds admits a symplectic
		structure, even without any invariance assumption.
		
		As a consequence, we show that in dimension four the underlying manifolds are precisely $L(d,1)\times S^1$, where $L(d,1)$ is a three-dimensional lens space. In dimension eight, we obtain families arising from seven-dimensional lens spaces and Aloff--Wallach spaces, including infinitely many distinct topological types.
	\end{abstract}

\section{Introduction}\label{section1}

A \emph{locally conformally symplectic} (LCS) structure on a manifold
$M$ is a pair $(\eta,\omega)$, where $\eta$ is a closed $1$-form,
called the \emph{Lee form}, and $\omega$ is a nondegenerate $2$-form
satisfying
\[
d\omega=\eta\wedge\omega.
\]
Locally, $\eta=df$, and the conformally related form $e^{-f}\omega$
is symplectic. Thus LCS geometry extends symplectic geometry.
An LCS structure is called \emph{exact} if there exists a $1$-form
$\lambda$ such that
\[
\omega=d\lambda-\eta\wedge\lambda.
\]

Locally conformally symplectic structures were introduced by Lee
\cite{Lee} and subsequently investigated by Vaisman \cite{Vaisman},
Banyaga \cite{banyaga1,banyaga2,banayaga3}, and many other authors
\cite{angella,BazzoniMarrero,bande,guedira,Le,marrero}.
They arise naturally in Hamiltonian mechanics as a generalization
of the symplectic framework \cite{Vaisman} and as the even-dimensional
transitive leaves of Jacobi manifolds \cite{guedira}.

Let $G$ be a connected Lie group, $K$ a closed subgroup, and
$\pi:G\to G/K$ the canonical projection. Suppose that $G/K$ carries
a $G$-invariant LCS structure $(\eta,\omega)$. Denote by
$\mathfrak g$ and $\mathfrak k$ the Lie algebras of $G$ and $K$,
respectively. The pullbacks $\pi^*\eta$ and $\pi^*\omega$ are
left-invariant forms on $G$. Their values at the identity define
\[
(\theta,\Omega)\in\mathfrak g^*\times\Lambda^2\mathfrak g^*.
\]
These infinitesimal data satisfy the following conditions:
\begin{enumerate}
	\item The LCS equations
	\begin{equation}\label{lcs-Lie}
		d_{\mathfrak g}\theta=0,\qquad
		d_{\mathfrak g}\Omega=\theta\wedge\Omega,
	\end{equation}
	where $d_{\mathfrak g}$ denotes the Chevalley--Eilenberg differential.
	
	\item The isotropy invariance conditions
	\begin{equation}\label{lcs-AdK}
		\operatorname{Ad}_k^*\theta=\theta,\qquad
		\operatorname{Ad}_k^*\Omega=\Omega,\qquad k\in K.
	\end{equation}
	
	\item The horizontality and nondegeneracy conditions
	\begin{equation}\label{lcs-kernel}
		\theta|_{\mathfrak k}=0,\qquad
		\ker\Omega=\mathfrak k,
	\end{equation}
	where $\ker\Omega=\{X\in\mathfrak g:\iota_X\Omega=0\}$.
\end{enumerate}
Conversely, every pair $(\theta,\Omega)$ satisfying
\eqref{lcs-Lie}--\eqref{lcs-kernel} determines a unique $G$-invariant
LCS structure on $G/K$ whose pullbacks have the prescribed values at
the identity.

In this paper, we classify compact connected homogeneous LCS
manifolds of dimension at least four with a nonzero Lee form.

The geometry of compact homogeneous locally conformal
symplectic manifolds was investigated by Vaisman
\cite{Vaisman}, who established important connections
with contact geometry and principal torus bundles over
homogeneous symplectic manifolds. Related results on
compact homogeneous contact manifolds and their
symplectic bases were obtained by Diaz Miranda and
Reventos \cite{Diaz}. More recently, Bazzoni and
Marrero \cite{BazzoniMarrero} established, under suitable
topological assumptions, a structure theorem for LCS
manifolds of the first kind, describing compact
examples as mapping tori of contact manifolds by
strict contactomorphisms. They also investigated
left-invariant LCS structures on Lie groups.

Our method is inspired by the approach of Zwart and
Boothby \cite{ZwartBoothby} to compact homogeneous
symplectic manifolds. We analyze the nilradical and
solvable radical of the \emph{original} transitive
Lie algebra. This leads to three explicit families
of pairs $(G,K)$, with $G$ connected and simply
connected and $K$ closed, for which $G/K$ is compact
and carries a $G$-invariant LCS structure with an
almost-effective $G$-action. We prove that these
families exhaust all possible homogeneous
presentations under our hypotheses (see Theorems \ref{main1}-\ref{main}).

Our second main result (Theorem \ref{thm:canonical-reduction}) provides an explicit global
geometric description of compact homogeneous LCS
manifolds with nonzero Lee form. We prove that every
such LCS manifold is LCS-isomorphic to the product
of a compact homogeneous contact manifold with a
circle, endowed with an explicitly determined LCS
structure. In particular, every such structure
admits a transitive, almost-effective action of
a compact Lie group of the form $S\times S^1$,
where $S$ is compact, simply connected, and
semisimple, with respect to which the LCS structure
is strongly exact, with an invariant primitive.

We also identify the contact manifold as a
principal circle bundle over a compact coadjoint
orbit. This construction is closely related to
the classical Boothby--Wang theorem \cite{Wang}.
More precisely, the contact form $\alpha_\beta$
is a connection form whose curvature satisfies
\[
d\alpha_\beta=-p^*\omega_{\mathrm{KKS}}.
\]
We also determine the period group of the
Kirillov--Kostant--Souriau symplectic form \cite{Kostant1970}
and show that its normalized cohomology class
is integral and primitive.

An interesting topological consequence of our
classification is that none of these compact
homogeneous LCS manifolds admits a symplectic
structure, even without any invariance assumption (See Proposition \ref{thm:topological-properties}).

While compact LCS manifolds without symplectic
structures are well known in the literature,
our result shows that this phenomenon holds
uniformly for the entire class of compact
homogeneous LCS manifolds with nonzero Lee form
considered here.

We now state our main results. Let $S$ be a compact, connected,
simply connected semisimple Lie group with Lie algebra $\mathfrak s$.
Choose $\beta\in\mathfrak s^*$ whose restriction to each simple ideal
is nonzero, and let
\[
H=S_\beta=\{h\in S:\operatorname{Ad}_h^*\beta=\beta\},
\qquad \mathfrak h=\operatorname{Lie}(H).
\]
Using an $\operatorname{Ad}(S)$-invariant inner product, write
$\beta=\langle b,\cdot\rangle$ for $b\in\mathfrak s$. Then
$H=\{ g\in S,\; \Ad_gb=b\}$ is  compact connected as the centralizer of a torus, $b\in\mathfrak h$, and
$\beta(b)=\langle b,b\rangle>0$. In particular,
$\beta|_{\mathfrak h}\ne0$.

Assume that $\beta|_{\mathfrak h}$ integrates to a surjective
character with connected kernel,
\begin{equation*}
	\phi_\beta:H\longrightarrow\mathbb T_\ell
	:=\mathbb R/\ell\mathbb Z,\qquad
	(d\phi_\beta)_e=\beta|_{\mathfrak h},\qquad\ell>0.
\end{equation*}
This integrability condition plays a central role in our
classification. In Section~\ref{section6}, we provide
a systematic method for constructing admissible covectors
$\beta$ satisfying this condition on compact, simply
connected semisimple Lie groups.

For an integer $d\ge1$, set
\[
F_d=\bigl\{[j\ell/d]_\ell:0\le j<d\bigr\},
\qquad H_d=\phi_\beta^{-1}(F_d).
\]
Then $H_d$ is a closed normal subgroup of $H$, and
$\operatorname{Lie}(H_d)=\mathfrak h\cap\ker\beta$.
The left invariant 1-form $\be^L$ associated to $\beta$ is $\operatorname{Ad}_{H_d}$-invariant and
vanishes on $\operatorname{Lie}(H_d)$; it therefore induces an
$S$-invariant contact $1$-form $\alpha_\beta$ on $S/H_d$.

Choose $a>0$ and $\alpha\in\mathbb R$, with $\alpha$ understood
modulo $\ell/d$. For the third family, let $V$ be a nonzero real
vector space and choose $D\in\operatorname{End}(V)$ and
$c\in V^*\setminus\{0\}$ satisfying
\begin{equation*}\label{eq:explicit-B-admissibility}
	D\ne0,\qquad
	\operatorname{span}\{(D^*)^jc:j\ge0\}=V^*,\qquad
	e^{aD}=\mathrm \mathrm{I}_V.
\end{equation*}
We have  $\{t\in\R\; \exp(tD)=\mathrm{I}_V \}=\tau\mathbb Z$ for a uniquely
determined $\tau>0$, and $a$ is a positive integer multiple of $\tau$.

The following three results give the homogeneous presentations
explicitly, including their invariant forms and their common
contact--circle realization.

\begin{theorem}[Model \textup{(A1)}: one-dimensional central radical]
	\label{main1}\label{prop:explicit-A1}
	Let $G_1=S\times\mathbb R$, and define
	\[
	K_1=\bigl\{(h,na):h\in H,\ n\in\mathbb Z,\
	\phi_\beta(h)-[n\alpha]_\ell\in F_d\bigr\}.
	\]
	On $\mathfrak g_1=\mathfrak s\oplus\mathbb RT$ let $t\in\G_1^*$ given by $t_{|\s}=0$ and $t(T)=1$. Put\footnote{Throughout the three models, $\beta$ and
		$d_{\mathfrak s}\beta$ are identified with their trivial
		extensions to $\mathfrak g$, that is, their pullbacks
		under the natural projection
		$\mathfrak g\to\mathfrak s$.}
	\begin{equation}\label{eq:explicit-A1-forms}
		\theta_1=t\esp 
		\Omega_1=d_{\mathfrak s}\beta+\beta\wedge t.
	\end{equation}
	Then $K_1$ is closed, $G_1/K_1$ is compact, and
	$(\theta_1,\Omega_1)$ induces a $G_1$-invariant LCS structure
	$(\eta_1,\omega_1)$ on $G_1/K_1$. The action is almost effective,
	and there is an LCS isomorphism
	\begin{equation}\label{eq:explicit-A1-reduction}
		(G_1/K_1,\eta_1,\omega_1)
		\simeq\bigl((S/H_d)\times S^1,\ a\,du,\
		d\alpha_\beta+a\,\alpha_\beta\wedge du\bigr).
	\end{equation}
\end{theorem}

\begin{theorem}[Model \textup{(A2)}: two-dimensional central radical]
	\label{main2}\label{prop:explicit-A2}
	Let $G_2=S\times\mathbb R_T\times\mathbb R_Z$, and define
	\[
	K_2=\bigl\{(h,na,z):h\in H,\ n\in\mathbb Z,\ z\in\mathbb R,\
	\phi_\beta(h)+[z-n\alpha]_\ell\in F_d\bigr\}.
	\]
	On $\mathfrak g_2=\mathfrak s\oplus\mathbb RT\oplus RZ$ let $(t,z)\in\G_2^*\times\G_2^*$ given by $t_{|\s\oplus \R Z}=0$,  $t(T)=1$, $z_{|\s\oplus \R T}=0$ and $z(Z)=1$. Put
	\begin{equation}\label{eq:explicit-A2-forms}
		\theta_2=t\esp
		\Omega_2=d_{\mathfrak s}\beta+(\beta+z)\wedge t.
	\end{equation}
	Then $K_2$ is closed, $G_2/K_2$ is compact, and these forms induce
	a $G_2$-invariant LCS structure $(\eta_2,\omega_2)$ with almost-effective
	action. Moreover,
	\begin{equation}\label{eq:explicit-A2-reduction}
		(G_2/K_2,\eta_2,\omega_2)
		\simeq\bigl((S/H_d)\times S^1,\ a\,du,\
		d\alpha_\beta+a\,\alpha_\beta\wedge du\bigr).
	\end{equation}
\end{theorem}

\begin{theorem}[Model \textup{(B)}: almost-abelian radical]
	\label{main3}\label{prop:explicit-B}
	Let
	\[
	G_B=S\times(\mathbb R_T\ltimes_D V),\qquad
	(t,v)(t',v')=(t+t',v+e^{tD}v'),
	\]
	and define
	\[
	K_B=\bigl\{(h,na,v):h\in H,\ n\in\mathbb Z,\ v\in V,\
	\phi_\beta(h)+[c(v)-n\alpha]_\ell\in F_d\bigr\}.
	\]
	On $\mathfrak g_B=\mathfrak s\oplus(\mathbb RT\ltimes_D V)$, let $t\in\G_B^*$ given by $t_{|\s\oplus V}=0$ and $t(T)=1$. Put
	\begin{equation}\label{eq:explicit-B-forms}
		\theta_B=t\esp
		\Omega_B=d_{\mathfrak s}\beta+(\beta+c)\wedge t.
	\end{equation}
	Then $K_B$ is closed, $G_B/K_B$ is compact, and the infinitesimal
	forms induce a $G_B$-invariant LCS structure $(\eta_B,\omega_B)$
	with almost-effective action. Furthermore,
	\begin{equation}\label{eq:explicit-B-reduction}
		(G_B/K_B,\eta_B,\omega_B)
		\simeq\bigl((S/H_d)\times S^1,\ a\,du,\
		d\alpha_\beta+a\,\alpha_\beta\wedge du\bigr),
	\end{equation}
	where $a\in\tau\mathbb Z_{>0}$.
\end{theorem}

These three families exhaust the compact homogeneous LCS
presentations considered here.

\begin{theorem}[Classification of compact homogeneous LCS structures]
	\label{main}\label{thm:explicit-compact-completeness}
	Let $M=G/K$ be a compact connected homogeneous manifold of dimension
	at least four, where $G$ is a connected, simply connected Lie group
	acting almost effectively. Suppose that $M$ carries a $G$-invariant
	LCS structure $(\eta,\omega)$ with nonzero Lee form.
	
	Then, up to an isomorphism of homogeneous LCS spaces,
	$(G/K,\eta,\omega)$ is one of the models described in
	Theorems~\ref{main1}--\ref{main3}, with admissible parameters.

\end{theorem}

\begin{theorem}[Canonical reduction and coadjoint-orbit geometry]
	\label{thm:canonical-reduction}
	Let $(M,\eta,\omega)$ be a compact homogeneous LCS manifold satisfying
	the assumptions of Theorem~\ref{main}. Then:
	\begin{enumerate}\item  There exist admissible
	data $(S,\beta,d,a)$ such that $(M,\eta,\omega)$ is LCS-isomorphic to
	\[
	\left((S/H_d)\times S^1,\ a\,du,\
	d\alpha_\beta+a\,\alpha_\beta\wedge du\right).
	\]
	Here $\alpha_\beta$ is the $S$-invariant contact form on $S/H_d$
	induced by the left-invariant one-form $\beta^L$ on $S$.
	The  LCS structure is preserved by the natural transitive,
	almost-effective action of $S\times S^1$ and is  exact,
	with $S\times S^1$-invariant primitive $\alpha_\beta$.
	
	\item  The natural projection
	$
	p:S/H_d\longrightarrow S/H\simeq\mathcal O_\beta
	$
	is a principal circle bundle, and
	$
	d\alpha_\beta=-p^*\omega_{\mathrm{KKS}},
	$
	where $\omega_{\mathrm{KKS}}$ is the Kirillov--Kostant--Souriau
	symplectic form on the coadjoint orbit $\mathcal O_\beta$.
	
	\item The  group of periods of $\omega_{\mathrm{KKS}}$ is precisely $\ell\mathbb Z$:
	\[
	\left\{\int_C\omega_{\mathrm{KKS}}:
	[C]\in H_2(\mathcal O_\beta;\mathbb Z)\right\}
	=\ell\mathbb Z.
	\]
	In particular, the real cohomology class
	$[\omega_{\mathrm{KKS}}]/\ell$ is the image of a unique
	primitive integral class
	$
	u_\beta\in H^2(\mathcal O_\beta;\mathbb Z).
	$
\end{enumerate}
\end{theorem}

The classification has explicit topological restrictions on $M$ (see Proposition \ref{thm:topological-properties}).
In dimension four, the underlying manifolds are precisely
$L(d,1)\times S^1$, where $L(d,1)$ is a three-dimensional lens space.
In dimension eight, the construction yields, among others,
products $W_{k,l}\times S^1$ involving Aloff--Wallach spaces;
these include infinitely many distinct topological types.

\paragraph{Organization of the paper.} The paper is organized as follows. Section \ref{section2} recalls the basic notions and results on locally conformally symplectic structures and homogeneous spaces. In Section~\ref{section3}, we study the infinitesimal structure of compact homogeneous LCS manifolds and obtain three models for the corresponding Lie algebra data. Section~\ref{section4} is devoted to the global integration of these models, with particular attention to the determination of the closed isotropy subgroups. In Section~\ref{section5}, we prove Theorems \ref{main1}-\ref{thm:canonical-reduction} and deduce some topological restrictions on compact homogeneous LCS manifolds (see Proposition \ref{thm:topological-properties}). 
Section~\ref{section6} provides a systematic
construction of the admissible covectors $\beta$
appearing in our classification. We establish
an explicit integrability criterion in terms
of the integral lattice of a maximal torus
and develop an algorithm for constructing
such covectors.
  Finally, Section~\ref{section7} provides explicit examples, including the complete description in dimension four and families of eight-dimensional homogeneous LCS manifolds arising from Aloff--Wallach spaces.

\section{LCS Structures and Preliminary  lemmas}\label{section2}

An LCS structure $(\omega,\eta)$ on $M$ is
called \emph{homogeneous} if there exists
a connected Lie group $G$ acting transitively
on $M$ and preserving both forms, that is,
\[
g^*\omega=\omega,\qquad
g^*\eta=\eta,\qquad g\in G.
\]
Fix a point $o\in M$ and denote by $K$
its isotropy subgroup. Then $K$ is a closed
subgroup of $G$, and the map
\[
G/K\longrightarrow M,\qquad
gK\longmapsto g\cdot o,
\]
is a $G$-equivariant diffeomorphism.
We may therefore identify $M$ with $G/K$
and regard $(\omega,\eta)$ as a
$G$-invariant LCS structure on this
homogeneous space.

We call a homogeneous LCS structure $(\eta,\om)$ \emph{exact} (resp.  \emph{strongly exact relative to $G$}) if there is a  1-form   $\mu$ (resp. $G$-invariant 1-form $\mu$) on $M$ such that
\begin{equation}\label{eq:strong}
	\omega=d\mu-\eta\wedge\mu.
\end{equation}
Every homogeneous LCS structure with nonvanishing Lee form is exact \cite[Proposition~2.4]{Vaisman}, but need not be strongly exact with respect to the given transitive group action. However, when $G$ is compact, averaging an ordinary primitive over $G$ yields a $G$-invariant primitive, so the structure is strongly exact.

Throughout this paper, we consider
homogeneous LCS structures on compact
homogeneous spaces $M=G/K$, where $G$
is connected and $K$ is closed.
We shall assume that the action of $G$
is \emph{almost effective}, meaning that
its kernel is discrete. Equivalently,
the isotropy Lie algebra $\mathfrak k=\operatorname{Lie}(K)$
contains no nonzero ideal of
$\mathfrak g=\operatorname{Lie}(G)$.

Let
$
\pi:G\longrightarrow G/K
$
be the canonical projection, and let
$(\omega,\eta)$ be a $G$-invariant LCS
structure on $G/K$.

The pullbacks $\pi^*\omega$ and $\pi^*\eta$
are left-invariant differential forms on
$G$. Identifying left-invariant forms
with their values at the identity, we
obtain a pair
$
(\Omega,\theta)\in
\Lambda^2\mathfrak g^*\times\mathfrak g^*.
$ This pair satisfies \eqref{lcs-Lie}--\eqref{lcs-kernel} and conversely, every pair $(\theta,\Omega)$ satisfying
\eqref{lcs-Lie}--\eqref{lcs-kernel} determines a unique $G$-invariant
LCS structure on $G/K$ whose pullbacks have the prescribed values at
the identity.

Consequently, the study of homogeneous
LCS structures on $G/K$ reduces to the
study of pairs $(\Omega,\theta)$
satisfying the above algebraic
and isotropy-invariance conditions.

\begin{remark}\label{rem:connected-isotropy}
	The $\operatorname{Ad}(K)$-invariance conditions
	\eqref{lcs-AdK} are automatically satisfied when
	$K$ is connected.
	
	Indeed, let $(\Omega,\theta)$ satisfy
	\eqref{lcs-Lie} and \eqref{lcs-kernel}.
	For any $X\in\mathfrak k$, Cartan's formula gives
	\[
	\begin{aligned}
		\mathcal L_X\Omega
		&=d_{\mathfrak g}(\iota_X\Omega)
		+\iota_X(d_{\mathfrak g}\Omega)
		=\iota_X(\theta\wedge\Omega)
		=\theta(X)\Omega-\theta\wedge\iota_X\Omega
		=0.
	\end{aligned}
	\]
	Similarly,
	\[
	\mathcal L_X\theta
	=d_{\mathfrak g}(\theta(X))
	+\iota_X(d_{\mathfrak g}\theta)=0.
	\]
	Thus both $\Omega$ and $\theta$ are infinitesimally
	invariant under $\mathfrak k$.
	
	If $K$ is connected, infinitesimal invariance
	implies invariance under the full adjoint action
	of $K$. Consequently, in this case,
	\eqref{lcs-AdK} follows from
	\eqref{lcs-Lie} and \eqref{lcs-kernel}.
\end{remark}

The following lemma strengthens
\cite[Lemma~3.3]{ZwartBoothby} and admits
a short elementary proof.  Its proof shows that neither the normality nor the connectedness assumptions in that lemma are needed.
It will play
a central role in our study.

\begin{lemma}[Compact-orbit lemma for unipotent representations]\label{Le:Compact}
	\label{lem:compact-unipotent}
	Let $G$ be a group and
	$\rho:G\to\mathrm{GL}(V)$ a finite-dimensional
	real representation. Let $N$ be a subgroup
	of $G$ such that $\rho(N)$ consists of
	unipotent transformations.
	
	If the orbit $\rho(G)v$ of some $v\in V$
	is bounded, then $v$ is fixed by $N$.
\end{lemma}

\begin{proof}
	Fix $a\in N$ and write
	$\rho(a)=I+U$, where $U^m=0$ for some
	positive integer $m$. For every
	integer $k\geq0$, we have
	\[
	\rho(a^k)v
	=(I+U)^kv
	=\sum_{j=0}^{m-1}
	\binom{k}{j}U^jv.
	\]
	The right-hand side is $P(k)$,
	the restriction of a polynomial map
	$P:\mathbb R\to V$ given by
	\[ P(t)=\sum_{j=0}^{m-1}\frac{t(t-1)\ldots(t-j+1)}{j!}U^jv. \]
	The sequence $(P(k))_{k\geq0}$ is bounded,
	since $P(k)\in\rho(G)v$. As a polynomial
	map bounded on the nonnegative integers
	must be constant, we obtain
	\[
	v=P(0)=P(1)=\rho(a)v.
	\]
	The proof is complete.
\end{proof}

The following lemma extends the infinitesimal
invariance assertion of
\cite[Lemma~4.1]{ZwartBoothby}
from bilinear forms to arbitrary alternating
multilinear forms.

\begin{lemma}\label{thm:p-Boothby}
	Let $G$ be a connected Lie group, $K$ a closed
	subgroup such that $G/K$ is compact, and
	$\nil$ the nilradical of $\g$.
	Let $p\geq 1$, let
	$\alpha\in\Lambda^p\g^*$ be
	$\Ad(K)$-invariant, and let
	$\mathfrak a_1,\ldots,\mathfrak a_p$
	be ideals of $\g$.
	Then, for every $Z\in\nil$ and
	$X_i\in\mathfrak a_i$, we have
	\begin{equation}\label{eq:p-Boothby}
		\sum_{i=1}^p
		\alpha(X_1,\ldots,[Z,X_i],\ldots,X_p)=0.
	\end{equation}
		In particular, taking
	$\mathfrak a_1=\cdots=\mathfrak a_p=\g$,
	we obtain $\ad_Z^*\al=0$ for any $Z\in\nil$.
\end{lemma}

\begin{proof}
	Since $G$ is connected, each ideal
	$\mathfrak a_i$ is preserved by
	$\Ad(G)$. Hence $G$ acts naturally on
	\[
	V=\mathfrak a_1^*\otimes\cdots
	\otimes\mathfrak a_p^*.
	\]
	Let $\rho:G\to\GL(V)$ denote the
	corresponding representation, and put
	\[
	v=\alpha|_{\mathfrak a_1\times
		\cdots\times\mathfrak a_p}.
	\]
	Since $\alpha$ is $\Ad(K)$-invariant,
	the stabilizer of $v$ contains $K$.
	Consequently, the orbit $\rho(G)v$
	is compact, being the continuous
	image of $G/K$.
	
	Now fix $Z\in\nil$. Since $\nil$
	is a nilpotent ideal of $\g$, we have
	$
	\ad_Z(\g)\subset\nil,
	$
	and the restriction of $\ad_Z$ to
	$\nil$ is nilpotent. It follows that
	$\ad_Z$ is nilpotent on $\g$.
	
	Consequently, its induced action
	$A=\rho_*(Z)$ on $V$ is also nilpotent.
	Thus
	\[
	\rho(\exp(tZ))=\exp(tA)
	\]
	is unipotent for every $t\in\mathbb R$.
	
	Applying Lemma~\ref{Le:Compact} to the
	one-parameter subgroup
	$\exp(\mathbb R Z)$, whose orbit of
	$v$ is contained in the compact set
	$\rho(G)v$, we obtain
	\[
	\rho_*(Z)v=0.
	\]
		By the definition of the induced
	representation, this identity reads
	\[
	\sum_{i=1}^p
	\alpha(X_1,\ldots,[Z,X_i],\ldots,X_p)=0
	\]
	for all $X_i\in\mathfrak a_i$.
	This proves the lemma.
\end{proof}

\section{Infinitesimal structure of compact homogeneous LCS manifolds}\label{section3}

Throughout this section, let $G$ be a connected
Lie group and $K$ a closed subgroup such that
$M=G/K$ is compact and the action of $G$ is
almost effective. Let $F: \wi G\too G$
be the universal covering group of $G$; since $M$ can naturally be identified with $\wi G/\wi K$, where $\wi K = F^{-1}(K)$, there is no loss of generality
in supposing that the group $G$ is simply connected.

Finally, we assume that $M$ has
dimension at least four and carries a
$G$-invariant LCS structure $(\eta,\om)$
with nonzero Lee form. 

Let $(\Om,\theta)\in
\Lambda^2\G^*\times\G^*$
be the associated infinitesimal pair,
let $\kalg$ be the Lie algebra of $K$,
and let $\nil$ denote the nilradical of
$\G=\mathrm{Lie}(G)$.
Recall that $(\Om,\theta)$ satisfies
\eqref{lcs-Lie}--\eqref{lcs-kernel}.

Applying Lemma~\ref{thm:p-Boothby} to $\Om$,
we obtain
\begin{equation}\label{eq:n-inv}
	\Om([Z,X],Y)+\Om(X,[Z,Y])=0,
	\qquad Z\in\nil,\quad X,Y\in\G.
\end{equation}
In particular, the Chevalley--Eilenberg
formula gives
\[
\CE\Om(Z,X,Y)=\Om(Z,[X,Y]).
\]
Combining this identity with
$\CE\Om=\theta\wedge\Om$, we obtain
\begin{equation}\label{eq:key}
	\begin{aligned}
		\Om(Z,[X,Y])
		={}&\theta(Z)\Om(X,Y)
		+\theta(X)\Om(Y,Z)
		+\theta(Y)\Om(Z,X),
	\end{aligned}
\end{equation}
for all $Z\in\nil$ and $X,Y\in\G$.

Moreover, for $U,V\in\nil$ and $X\in\G$,
equation \eqref{eq:n-inv} implies
\[
\Om([U,V],X)=\Om(U,[V,X]).
\]
Applying \eqref{eq:key}, we therefore get
\begin{equation}\label{eq:comm}
	\begin{aligned}
		\Om([U,V],X)
		={}&\theta(U)\Om(V,X)
		-\theta(V)\Om(U,X)
		+\Om(U,V)\theta(X).
	\end{aligned}
\end{equation}
These identities are the key ingredients
in the following structural results.

\begin{theorem}[Abelianity of the nilradical]
	\label{thm:nilrad}
	Under the hypotheses of this section, 
	 the nilradical $\nil$ of $\G$
	is abelian:
	$
	[\nil,\nil]=0.
	$
	Furthermore, if $\theta|_{\nil}=0$ then
	$
	\Om(\nil,\nil)=0.
	$
\end{theorem}

\begin{proof}
	We distinguish two cases.
	
	\smallskip
	\noindent
	\textbf{Case 1: $\theta|_{\nil}=0$.}
	
	Since $\nil$ is an ideal of $\G$
	and $G$ is connected, $\nil$ is
	preserved by $\Ad(G)$. Consider the
	induced representation
	$
	\rho:G\longrightarrow\GL(\Lambda^2\nil^*)
	$
	defined by
	\[
	(\rho(g)\sigma)(U,V)
	=\sigma(\Ad_{g^{-1}}U,\Ad_{g^{-1}}V),\;\;g\in G,\;\sigma\in \Lambda^2\nil^*,\; U,V\in\nil.
	\]
	Set
	$
	\sigma=\Om|_{\nil\times\nil}.
	$
	Since $\theta$ vanishes on $\nil$,
	equation \eqref{eq:comm} reduces to
	\begin{equation}\label{eq:solvable-sigma-comm}
		\Om([U,V],X)=\sigma(U,V)\theta(X),
	\end{equation}
	for all $U,V\in\nil$ and $X\in\G$.
	By \eqref{eq:n-inv}, we have
	\[
	\Om([X,U],V)
	=-\sigma(U,V)\theta(X)\esp
	\Om(U,[X,V])
	=-\sigma(U,V)\theta(X).
	\]
	Consequently, the infinitesimal action
	of $\G$ on $\sigma$ satisfies
	\begin{equation*}\label{eq:solvable-semi-invariance}
		\begin{aligned}
			(\rho_*(X)\sigma)(U,V)
			&=-\sigma([X,U],V)
			-\sigma(U,[X,V])
			=2\theta(X)\sigma(U,V).
		\end{aligned}
	\end{equation*}
	Thus
	$
	\rho_*(X)\sigma=2\theta(X)\sigma,
	$
	and hence
	\begin{equation}\label{exp}
		\rho(\exp(tX))\sigma
		=e^{2t\theta(X)}\sigma.
	\end{equation}
		Since $\Om$ is $\Ad(K)$-invariant,
	so is $\sigma$. Therefore the stabilizer
	$G_\sigma$ contains $K$, and the orbit
	$
	\rho(G)\sigma\simeq G/G_\sigma
	$
	is compact, being the continuous image
	of the compact space $G/K$.
	
	As $\theta\neq0$, we can choose
	$X\in\G$ with $\theta(X)=1$.
	Equation \eqref{exp} then gives
	\[
	e^{2t}\sigma\in\rho(G)\sigma,
	\qquad t\in\mathbb R.
	\]
	The compactness of this orbit forces
	$\sigma=0$. Hence
	$
	\Om(\nil,\nil)=0.
	$
	It follows from
	\eqref{eq:solvable-sigma-comm} that
	$
	\Om([\nil,\nil],\G)=0.
	$
	Since $\ker\Om=\kalg$, we deduce that
	$
	[\nil,\nil]\subseteq\kalg.
	$
	But $[\nil,\nil]$ is an ideal of $\G$.
	By almost-effectivity,
	$\kalg$ contains no nonzero ideal of
	$\G$. Consequently,
	$
	[\nil,\nil]=0.
	$
	
	\smallskip
	\noindent
	\textbf{Case 2: $\theta|_{\nil}\neq0$.}
	
	Choose $V\in\nil$ such that
	$\theta(V)=1$, and put
	$
	\nil_0=\nil\cap\ker\theta.
	$
	Then
	$
	\nil=\mathbb RV\oplus\nil_0.
	$
	
	Since $\ker\Om=\kalg$ and
	$\theta|_{\kalg}=0$, the forms
	$\Om$ and $\theta$ induce a
	nondegenerate skew-symmetric form
	$\overline{\Om}$ and a linear form
	$\overline{\theta}$ on the vector space
	$
	\overline{\G}=\G/\kalg.
	$
	Let $\bar E\in\overline{\G}$
	be the unique vector satisfying
	\begin{equation}\label{eq:nil-E}
		\overline{\Om}(\bar E,\cdot)
		=\overline{\theta}.
	\end{equation}
	
	Since $V\in\nil$, the endomorphism
	$\ad_V$ is nilpotent because $[V,\G]\subset\nil$ and the restriction of $\ad_V$ to $\nil$ is nilpotent. Moreover,
	\eqref{eq:n-inv} shows that
	$\ad_V$ preserves $\ker\Om=\kalg$.
	It therefore induces a nilpotent
	endomorphism
	\[
	T:\overline{\G}\longrightarrow
	\overline{\G},
	\qquad
	T\bar X=\overline{[V,X]}.
	\]
	By \eqref{eq:n-inv} and the closedness
	of $\theta$, we have
	\begin{equation*}\label{eq:nil-T}
		\overline{\Om}(T\bar X,\bar Y)
		+\overline{\Om}(\bar X,T\bar Y)=0\esp 
		\overline{\theta}\circ T=0.
	\end{equation*}
	Taking $\bar X=\bar E$ in the first
	identity and using
	\eqref{eq:nil-E}, we obtain
	\[
	\overline{\Om}(T\bar E,\bar Y)
	=-\overline{\theta}(T\bar Y)=0.
	\]
	By nondegeneracy,
	\begin{equation}\label{eq:nil-TE}
		T\bar E=0.
	\end{equation}
	Now let $U\in\nil_0$. Since
	$\theta(U)=0$ and $\theta(V)=1$,
	equation \eqref{eq:comm} gives
	\[
	\Om([V,U],X)
	=\Om(U,X)+\Om(V,U)\theta(X).
	\]
	Passing to the quotient, we obtain
	$
		T\bar U=\bar U+\Om(V,U)\bar E.
	$
	Applying $T$ to this identity
	and using \eqref{eq:nil-TE}, we get
	$
	T^2\bar U=T\bar U.
	$
	Since $T$ is nilpotent, this implies
	$T\bar U=0$. Therefore
	\begin{equation}\label{eq:nil-U}
		[V,U]\in\kalg,
		\qquad
		\bar U=-\Om(V,U)\bar E.
	\end{equation}
	In particular, for any
	$U_1,U_2\in\nil_0$, their images
	in $\overline{\G}$ are collinear
	with $\bar E$. Hence
	\[
	\Om(U_1,U_2)
	=\overline{\Om}(\bar U_1,\bar U_2)
	=0.
	\]
	Applying \eqref{eq:comm} once again,
	and using $\theta(U_1)=\theta(U_2)=0$,
	we obtain
	\[
	\Om([U_1,U_2],X)=0,
	\qquad X\in\G.
	\]
	Thus
	$
	[\nil_0,\nil_0]\subseteq\kalg.
	$
	Together with \eqref{eq:nil-U},
	this yields
	\[
	[\nil,\nil]
	=[\nil_0,\nil_0]+[V,\nil_0]
	\subseteq\kalg.
	\]
	Finally, $[\nil,\nil]$ is an ideal
	of $\G$. By almost-effectivity,
	we conclude that
	$
	[\nil,\nil]=0.
	$
	This completes the proof.
\end{proof}

\begin{theorem}[Splitting of the Levi decomposition]
	\label{prop:mixed-Levi-splits}
	Under the hypotheses of this section, let
	$\mathfrak n$ and $\mathfrak r$ denote the
	nilradical and solvable radical of $\mathfrak g$,
	respectively. Then every Levi factor
	$\mathfrak s$ of $\mathfrak g$ centralizes
	$\mathfrak r$, that is,
	\[
	[\mathfrak s,\mathfrak r]=0.
	\]
	Consequently, the Levi factor is unique, and
	the Levi decomposition is a direct sum of ideals:
	\[
	\mathfrak g=\mathfrak s\oplus\mathfrak r.
	\]
\end{theorem}

\begin{proof}
	Fix a Levi factor $\mathfrak s$, so that
	$\mathfrak g=\mathfrak s\ltimes\mathfrak r$.
	We first prove that
	$
	[\mathfrak s,\mathfrak n]=0.
	$
	
	Put $\mathfrak n_0=\mathfrak n\cap\ker\theta$.
	Since $\theta$ is closed and $\mathfrak s$
	is perfect, we have
	$
	\theta(\mathfrak s)=0.
	$
	Thus, for any $Z\in\mathfrak n_0$ and
	$X,Y\in\mathfrak s$, equation~\eqref{eq:key}
	gives
	$
	\Omega(Z,[X,Y])=0.
	$
	Since $[\mathfrak s,\mathfrak s]=\mathfrak s$,
	we obtain
	\begin{equation}\label{eq:mixed-n0s}
		\Omega(\mathfrak n_0,\mathfrak s)=0.
	\end{equation}
	Now let $Z\in\mathfrak n$ and
	$X\in\mathfrak s$. Since $\mathfrak n$
	is an ideal and $\theta$ is closed,
	$
	[Z,X]\in\mathfrak n_0.
	$
	Hence, by \eqref{eq:mixed-n0s},
	$
	\Omega([Z,X],\mathfrak s)=0.
	$
	
	On the other hand, since $[\G,\mathfrak r]\subset\nil$ (see \cite[Corollary~5.4.15 and Remark~5.4.16]{HilgertNeeb}) and $\theta$ is closed, for any
	$Y\in\mathfrak r$, we have
	$[Z,Y]\in\mathfrak n_0$. Therefore,
	equations~\eqref{eq:n-inv} and
	\eqref{eq:mixed-n0s} yield
	\[
	\Omega([Z,X],Y)
	=-\Omega(X,[Z,Y])=0.
	\]
	Consequently,
	$
	\Omega([\mathfrak s,\mathfrak n],
	\mathfrak g)=0,
	$
	and hence
	$
		[\mathfrak s,\mathfrak n]
		\subseteq\ker\Omega=\mathfrak k.
	$
	
	We claim that
	$I=[\mathfrak s,\mathfrak n]$
	is an ideal of $\mathfrak g$.
	Indeed, Jacobi's identity and the fact
	that $\mathfrak n$ is an ideal give
	$
	[\mathfrak s,I]\subseteq I.
	$
	Moreover, for $Y\in\mathfrak r$,
	$X\in\mathfrak s$, and $Z\in\mathfrak n$,
	Jacobi's identity yields
	\[
	[Y,[X,Z]]
	=[[Y,X],Z]+[X,[Y,Z]].
	\]
	Since
	$
	[\mathfrak g,\mathfrak r]\subseteq\mathfrak n
	$
	and $\mathfrak n$ is abelian by
	Theorem~\ref{thm:nilrad}, the first term
	vanishes, while the second belongs to $I$.
	Thus
	$
	[\mathfrak r,I]\subseteq I,
	$
	proving that $I$ is an ideal of $\mathfrak g$
	and by  almost-effectivity,
	we conclude that
	$
	[\mathfrak s,\mathfrak n]=0.
	$
	
	Finally, since $\mathfrak s$ is perfect and
	$[\mathfrak s,\mathfrak r]\subseteq\mathfrak n$,
	Jacobi's identity gives
	\[
	\begin{aligned}
		[\mathfrak s,\mathfrak r]
		&=[[\mathfrak s,\mathfrak s],\mathfrak r]
		\subseteq
		[\mathfrak s,[\mathfrak s,\mathfrak r]]
		\subseteq[\mathfrak s,\mathfrak n]
		=0.
	\end{aligned}
	\]
	Therefore
	$
	\mathfrak g=\mathfrak s\oplus\mathfrak r
	$
	is a direct sum of ideals.
	The uniqueness of
	$\mathfrak s$ follows from the Levi--Malcev
	conjugacy theorem.
\end{proof}

We begin by determining the restriction of
$\Om$ to the Levi factor and its mixed component.

\begin{lemma}\label{lem:nf-beta}
	Under the preceding assumptions, let
	$\G=\mathfrak s\oplus\mathfrak r$
	be the Levi decomposition.
	There exists a unique $\be\in\mathfrak s^*$
	such that, for any $X,Y\in\mathfrak s$ and $U\in\mathfrak r$,
	\begin{equation}\label{eq:nf-beta}
		\Om(X,Y)=-\be([X,Y])
			\esp 
			\Om(X,U)=\be(X)\theta(U),
	\end{equation}
\end{lemma}

\begin{proof}
	Since $\mathfrak s$ is perfect and
	$\theta$ is closed, we have
	$\theta|_{\mathfrak s}=0$. By virtue of \eqref{lcs-Lie},
	 the restriction of $\Om$ to
	$\mathfrak s$ is a Chevalley--Eilenberg
	$2$-cocycle.
	By  the second Whitehead lemma, there exists
	$\be\in\mathfrak s^*$ such that
	\[
	\Om(X,Y)=-\be([X,Y]),
	\qquad X,Y\in\mathfrak s.
	\]
	The perfectness of $\mathfrak s$ ensures
	the uniqueness of $\be$.
	
	Now let $X,Y\in\mathfrak s$ and
	$U\in\mathfrak r$. Since
	$[\mathfrak s,\mathfrak r]=0$,
	the LCS equation $d\Om=\theta\wedge\Om$ gives
	\[
	-\Om([X,Y],U)
	=\theta(U)\Om(X,Y)
	=-\theta(U)\be([X,Y]).
	\]
	As $\mathfrak s=[\mathfrak s,\mathfrak s]$,
	we obtain
	$
	\Om(X,U)=\be(X)\theta(U),
	$
	as required.
\end{proof}

Let $S$ be the connected normal subgroup of $G$ associated to $\mathfrak s$ and $S_\beta=\{g\in S,\;\Ad_g^*\be=\be \}$ the isotropy subgroup of $\be$. The Lie subalgebra of $S_\beta$ is given by
\[
\mathfrak s_\be
=\{X\in\mathfrak s:
\be([X,\mathfrak s])=0\}.
\]
We denote by $\mathcal O_\beta=S/S_\beta$ the coadjoint orbit of $\be$.

\begin{lemma}\label{lem:nf-compact}
	Under the preceding assumptions,
	$\mathfrak s$ is compact semisimple,
	and $\be$ does not vanish on any
	simple ideal of $\mathfrak s$.
	Moreover, if $\mathfrak s\neq0$, then
	$
	\be|_{\mathfrak s_\be}\neq0.
	$
\end{lemma}

\begin{proof}
	We first observe that $\be$ cannot vanish
	on any simple ideal $\mathfrak a$ of
	$\mathfrak s$. Indeed, otherwise
	\eqref{eq:nf-beta} would imply
	$
	\Om(\mathfrak a,\G)=0,
	$
	so that the nonzero ideal $\mathfrak a$
	would be contained in $\kalg=\ker\Om$,
	contradicting almost-effectivity.
	
	Next, since $\Om$ is $\Ad(K)$-invariant
	and $\be$ is uniquely determined by
	\eqref{eq:nf-beta}, the linear form
	$\be$ is fixed by the induced coadjoint
	action of $K$. Consequently, its
	coadjoint orbit $\mathcal O_\be$ is
	compact, being the continuous image
	of the compact homogeneous space $G/K$.
	
	Suppose that $\mathfrak s$ contains
	a noncompact simple ideal $\mathfrak s_0$,
	and put
	$
	\be_0=\be|_{\mathfrak s_0}.
	$
	
	Let
	$
	r:\mathfrak s^*\longrightarrow\mathfrak s_0^*,
	\;
	r(\alpha)=\alpha|_{\mathfrak s_0},
	$
	be the restriction map.
	Since $\mathfrak s_0$ is an ideal,
	the coadjoint orbit of
	$\be_0=r(\be)$ is contained in
	$r(\mathcal O_\be)$.
	As $\mathcal O_\be$ is compact,
	$r(\mathcal O_\be)$ is compact,
	and hence the coadjoint orbit of
	$\be_0$ is bounded.
	
	By \cite[Theorem~2]{Rawnsley}, every
	bounded coadjoint orbit of a noncompact
	simple real Lie algebra is reduced
	to the origin. Thus $\be_0=0$,
	contradicting our initial observation.
	Therefore $\mathfrak s$ is compact
	semisimple.
	
	Finally, suppose $\mathfrak s\neq0$.
	Let $B$ be the Killing form of
	$\mathfrak s$, and let $b\in\mathfrak s$
	be the unique element satisfying
	$
	\be(X)=B(b,X)
	$ for any $ X\in\mathfrak s$.
	Since $\mathfrak s$ is compact,
	$-B$ is positive definite. Moreover,
	\[
	\mathfrak s_\be
	=\{X\in\mathfrak s:[X,b]=0\}.
	\]
	In particular, $b\in\mathfrak s_\be$.
	Since $\be\neq0$, we have $b\neq0$,
	and therefore
	\[
	\be(b)=B(b,b)<0.
	\]
	Hence
	$
	\be|_{\mathfrak s_\be}\neq0,
	$
	which completes the proof.
\end{proof}

\begin{theorem}[Infinitesimal normal forms]
	\label{thm:nf-main}
	Under the preceding assumptions, let
	$\be\in\mathfrak s^*$ be the unique
	linear form determined by
	\eqref{eq:nf-beta}. Then $\mathfrak s$ is
	compact semisimple, $\mathfrak{s}\not=0$ and $\be$ is
	nonzero on every simple ideal and
	precisely one of the following
	three models occurs, where
	$\mathfrak k=\ker\Om$.
	
	\smallskip
	\noindent
	\textup{\textbf{(A1) One-dimensional central radical.}}
	There exists $T\in Z(\G)$ such that
	\begin{equation}\label{eq:nf-A1}
		\begin{cases}
			\G=\mathfrak s\oplus\mathbb RT,\;\kalg=\mathfrak s_\be\cap\ker\be
			,\\
			\Om_{|\mathfrak s\times \mathfrak s}=d_{\mathfrak s}\be, \; \Om(\mathfrak s,T)=\be(\mathfrak s),\;
			\theta_{|\mathfrak s}=0,\; \theta(T)=1.
		\end{cases}
	\end{equation}

	\smallskip
	\noindent
	\textup{\textbf{(A2) Two-dimensional central radical.}}
	There exists a central basis $(T,Z)$ of $Z(\G)$,
	such that
	\begin{equation}\label{eq:nf-A2}
		\begin{cases}
			\G=\mathfrak s\oplus
			\mathbb RT\oplus\mathbb RZ,\; \kalg=
			\{X-\be(X)Z:X\in\mathfrak s_\be\},\\
			
			\Om_{|\mathfrak s\times \mathfrak s}=d_{\mathfrak s}\be,\;
			\Om(\mathfrak s,T)=\be(\mathfrak s),\; \Om(\mathfrak s,Z)=0,\; \Om(Z,T)=1,\;  \theta_{|\mathfrak s\oplus\R Z}=0,\;\theta(T)=1.
			
		\end{cases}
	\end{equation}
	
	\smallskip
	\noindent
	\textup{\textbf{(B) Almost-abelian radical.}}
	There exist a nonzero vector space
	$\nil$, a nonzero endomorphism
	$D\in\operatorname{End}(\nil)$,
	and a nonzero $c\in\nil^*$ such that
	\begin{equation}\label{eq:nf-B}
		\begin{cases}
			\G=\mathfrak s\oplus
			(\mathbb RT\ltimes_D\nil),\;
			[T,Z]=DZ,\; [\nil,\nil]=0,\;Z\in\nil,\\ \kalg=
			\{X+Z:X\in\mathfrak s_\be,\ Z\in\nil,\
			\be(X)+c(Z)=0\},\; \kalg\cap\nil=\ker c,\\
			\Om_{|\mathfrak s\times \mathfrak s}=d_{\mathfrak s}\be,\;
			\Om(\mathfrak s,T)=\be(\mathfrak s),\;\Om(\mathfrak s,\nil)=0,\;
			\Om(\nil,\nil)=0,\;\Om(\nil,T)=c(\nil),\\
			\theta_{|\mathfrak s\oplus\nil}=0,\; \theta(T)=1.
		\end{cases}
	\end{equation}
	Moreover, $c$ is cyclic for the dual endomorphism $D^*$, i.e., $\Span\{(D^*)^nc,n\geq0 \}=\nil^*$ and the period group of $D$ is a
	nontrivial discrete subgroup of
	$\mathbb R$. More precisely, there exists $\tau>0$ such that
	\[
	P_D=
	\{t\in\mathbb R:
	\exp(tD)=\operatorname{Id}_{\mathfrak n}\}
	=\tau\mathbb Z.
	\]

	Finally, in  all the three model, the dimension of the
	homogeneous space satisfies
	$
		\dim(G/K)=\dim\mathcal O_\be+2.
	$
\end{theorem}

\begin{proof}
	By Theorem~\ref{prop:mixed-Levi-splits},
	$
	\G=\mathfrak s\oplus\mathfrak r,
	\; [\mathfrak s,\mathfrak r]=0.
	$
	By Theorem~\ref{thm:nilrad}, the
	nilradical $\nil$ is abelian.
	Lemma~\ref{lem:nf-beta} determines
	the restriction of $\Om$ to
	$\mathfrak s$ and its mixed component,
	while Lemma~\ref{lem:nf-compact}
	establishes the assertions concerning
	$\mathfrak s$ and $\be$. Moreover, $\theta([\G,\G])=0$ and, in particular, since $\mathfrak{s}$ is perfect $\theta_{|\mathfrak{s}}=0$. We will also use the equations \eqref{eq:n-inv}-\eqref{eq:comm}.
	
	We first determine the radical and
	then the corresponding normal forms.
	
	\medskip
	\noindent
	\textbf{1. Structure of the radical.}
	
	Suppose first that
	$\theta|_{\nil}\neq0$, and put
	$\nil_0=\nil\cap\ker\theta$.
	Equation~\eqref{eq:key} and the
	perfectness of $\mathfrak s$ give
	$\Om(\nil_0,\mathfrak s)=0.$
	
	Let $W=[\mathfrak r,\nil]$.
	Since $\nil$ is abelian,
	\eqref{eq:n-inv} implies
	$
	\Om(W,\nil)=0.
	$
	Moreover, $W\subseteq\nil_0$, so
	$\Om(W,\mathfrak s)=0$.
	
	For $w\in W$ and
	$X,Y\in\mathfrak r$, we have
	$[X,Y]\in\nil$,  $\theta(w)=0$ and $\Om(w,[X,Y])=0$.
	Equation~\eqref{eq:key} therefore yields
	\[
	\theta(X)\Om(w,Y)
	=\theta(Y)\Om(w,X).
	\]
	But $\theta_{|\mathfrak r}\not=0$ and we can choose $Y\in \mathfrak r$ such that $\theta(Y)=1$. This implies that there exists $a_w\in\R$ such that, for any $X\in \mathfrak r$,
	$\Om(w,X)=a_w\theta(X)$. 
	
	We have also $\theta|_{\nil}\neq0$ and
	$\Om(w,\nil)=0$ so $a_w=0$. We have thus shown  that $\Om(W,\mathfrak s)=\Om(W,\mathfrak r)=0$ and
	consequently,
	$
	W\subseteq\ker\Om=\kalg.
	$
	As $W$ is an ideal of $\G$,
	almost-effectivity gives $W=0$.
	
	Since
	$
	[\mathfrak r,\mathfrak r]\subseteq\nil$
	and
	$[\mathfrak r,\nil]=0,$
	the radical $\mathfrak r$ is nilpotent.
	Hence $\mathfrak r=\nil$, and therefore
	$\mathfrak r=\nil=Z(\G).
	$
	
	Choose $T\in\nil$ with $\theta(T)=1$.
	The LCS equation $d\Om=\theta\wedge\Om$ applied to $(T,X,Y)$ with $X,Y\in\nil_0$ implies
	$
	\Om(\nil_0,\nil_0)=0.
	$
	We already know that
	$\Om(\nil_0,\mathfrak s)=0$.
	The  kernel of the linear map
	\[
	\kappa:\nil_0\longrightarrow\mathbb R,
	\qquad U\longmapsto\Om(U,T),
	\]is a
	central ideal contained in $\kalg$
	and must therefore be zero. 
	It follows that
	\begin{equation}\label{dim}
		1\leq\dim\mathfrak r=\dim\nil
		=1+\dim\nil_0\leq2.
	\end{equation}
	Finally, $\G=\mathfrak{s}\oplus\nil$ and $\dim\nil\leq 2$. It follows that $\s\not=0$ otherwise $\dim M\leq2$.
	
	Suppose now that $\theta|_{\nil}=0$.
	By Theorem~\ref{thm:nilrad},
	$
	\Om(\nil,\nil)=0.
	$
	Equation~\eqref{eq:key} and the perfectness of $\mathfrak s$ give
	$\Om(\nil,\mathfrak s)=0$.
	
	Choose $T\in\mathfrak r$ with
	$\theta(T)=1$ and define
	\[
	c(Z)=\Om(Z,T),\qquad Z\in\nil.
	\]
	Since
	$[\G,\mathfrak r]=[\mathfrak r,\mathfrak r]\subseteq\nil$, for any $Z\in\nil$ and $X\in\G$, $\Om(Z,[T,X])=0$ and
	equation~\eqref{eq:key} yields
	\begin{equation}\label{eq:nf-proof-c}
		\Om(Z,X)=c(Z)\theta(X),
		\qquad Z\in\nil,\ X\in\G.
	\end{equation}
	In particular, $
	\nil\cap\kalg=\ker c.$
	Almost-effectivity implies $c\neq0$.
	
	Set
	$\mathfrak h=\mathfrak r\cap\ker\theta$.
	Then $\mathfrak h$ is an ideal
	containing $\nil$.
	For $X\in\mathfrak h$ and
	$Z\in\nil$, the definition of $c$,  equations~\eqref{eq:n-inv} and \eqref{eq:nf-proof-c}
	give
	\[
	c([Z,X])=\Om([Z,X],T)=-\Om(X,[Z,T])=c([Z,T])\theta(X)=0.
	\]
	Hence
	$
	[\mathfrak h,\nil]\subseteq\kalg.
	$
	Since this bracket space is an
	ideal of $\G$, almost-effectivity
	implies
	$
	[\mathfrak h,\nil]=0.
	$
	As
	$[\mathfrak h,\mathfrak h]\subseteq\nil$,
	the ideal $\mathfrak h$ is nilpotent.
	Therefore $\mathfrak h=\nil$.
	
	We conclude that
	\begin{equation*}\label{eq:D}
		\mathfrak r=\mathbb RT\ltimes_D\nil,
		\qquad
		D=\operatorname{ad}_T|_{\nil}.
	\end{equation*}
	Moreover, $D\neq0$, since otherwise
	$\mathfrak r$ would be abelian and
	hence equal to $\nil$, contradicting
	$\theta(T)=1$. If $\s=0$, $\G=\mathbb RT\ltimes_D\nil$ and $\Om(\nil,\nil)=0$. Thus $\rank\Om\leq 2$ which contradicts $\dim M\geq4$.
	
	\medskip
	\noindent
	\textbf{2. Determination of the three models.}
	
	{\bf First case: $\theta|_{\nil}\neq0$}. By \eqref{dim},  we have $\theta_{|\mathfrak{s}}=0$ and $1\leq \dim \mathfrak r\leq 2$. We distinguish two cases.
	\begin{itemize}
		\item[$(a)$]	$\dim\mathfrak r=1$, choose
		$T\in \mathfrak r= Z(\G)$ with $\theta(T)=1$. Then $\G=\mathfrak{s}\oplus\R T$ and,
		by Lemma~\ref{lem:nf-beta},	
		\[\Om_{|\mathfrak s\times \mathfrak s}=d_{\mathfrak s}\be,\;\Om(\mathfrak s,T)=\beta(\mathfrak s).
		\]
		Here $\be\neq0$, since otherwise
		$\Om=0$, contradicting
		$\ker\Om\subseteq\ker\theta$.
		
		Thus $\mathfrak s\neq0$, and
		Lemma~\ref{lem:nf-compact} gives
		$\be|_{\mathfrak s_\be}\neq0$.
		For $X\in\mathfrak s$ and
		$a\in\mathbb R$, the condition
		$X+aT\in\ker\Om$ is therefore
		equivalent to
		\[
		a=0,\qquad
		X\in\mathfrak s_\be\cap\ker\be.
		\]
		This gives Model~\textup{(A1)}.
		
		\item[$(b)$]
		$\dim\mathfrak r=2$. Choose a
		central basis $(T,Z)$ of $\mathfrak{r}$ such that
		$
		\theta(T)=1,\;\theta(Z)=0.
		$  By Lemma~\ref{lem:nf-beta},
		\[ \Om_{|\mathfrak s\times\mathfrak s}=d_{\mathfrak s}\be,\;
		\Om(\mathfrak{s},T)=\beta(\mathfrak{s})\esp \Om(\mathfrak{s},Z)=0. \]
		Almost-effectivity forces
		$\Om(Z,T)\neq0$, since otherwise the central ideal $\R Z$ is contained in $\kalg$.
		After rescaling $Z$, we may assume
		$\Om(Z,T)=1$.
		A direct calculation gives
		\[
		\ker\Om
		=\{X-\be(X)Z:
		X\in\mathfrak s_\be\},
		\]
		which is Model~\textup{(A2)}.
		
	\end{itemize}
	
	{\bf Second case $\theta|_{\nil}=0$}. We have seen \eqref{eq:D} that
	\[ \mathfrak r=\mathbb RT\ltimes_D\nil,
	\qquad
	D=\operatorname{ad}_T|_{\nil}\esp D\not=0. \]
	From Lemma~\ref{lem:nf-beta}
	and \eqref{eq:nf-proof-c}, we obtain
	\[\Om_{|\mathfrak s\times \mathfrak s}=d_{\mathfrak s}\be,\;
	\Om(\mathfrak s,T)=\be(\mathfrak s),\;\Om(\mathfrak s,\nil)=0,\;
	\Om(\nil,\nil)=0,\;\Om(\nil,T)=c(\nil)\esp c\in\nil^*\setminus\{0\}.
	\]
	Since $c\neq0$, an element
	$X+aT+Z$, with
	$X\in\mathfrak s$ and $Z\in\nil$,
	belongs to $\ker\Om$ if and only if
	\[
	a=0,\qquad
	X\in\mathfrak s_\be,\qquad
	\be(X)+c(Z)=0.
	\]
	This gives the asserted expression
	for $\kalg$ in Model~\textup{(B)}.
	
	\medskip
	\noindent

	\medskip
	\noindent
	\textbf{3. Cyclicity of $c$.}
	
	Let $W\subseteq\mathfrak n$
	be a $D$-invariant subspace contained
	in $\ker c$.
	Since $\mathfrak n$ is abelian and
	$[\mathfrak s,\mathfrak n]=0$,
	$W$ is an ideal of $\mathfrak g$.
	Moreover,
	\[
	W\subseteq\ker c=\nil\cap\kalg\subseteq\mathfrak k.
	\]
	Almost-effectivity therefore implies
	$W=0$.
	
	Now consider
	$
	E=\operatorname{span}
	\{c,D^*c,(D^*)^2c,\ldots\}
	\subseteq\mathfrak n^*.
	$
	Its annihilator $E^\circ$ is
	$D$-invariant and contained
	in $\ker c$. Hence $E^\circ=0$,
	and therefore
	\begin{equation}\label{eq:mixed-cyclic}
		\mathfrak n^*=
		\operatorname{span}
		\{c,D^*c,(D^*)^2c,\ldots\}.
	\end{equation}
	
	\medskip
	\noindent
	\textbf{4. Periodicity of $D$.}
	
	Since $G$ is simply connected,
	the Lie algebra character
	$\theta$ integrates to a group
	homomorphism
	\[
	\chi:G\longrightarrow\mathbb R,
	\qquad d\chi_e=\theta.
	\]
	
	We claim that
	\begin{equation}\label{eq:period-ad}
		\Ad_g|_{\mathfrak n}
		=\exp(\chi(g)D),
		\qquad g\in G.
	\end{equation}
	Indeed, for
	$X=S+aT+Z\in\mathfrak g$,
	with $S\in\mathfrak s$ and
	$Z\in\mathfrak n$, we have
	$
	\ad_X|_{\mathfrak n}
	=aD=\theta(X)D.
	$
	Thus the two representations
	\[
	g\longmapsto\Ad_g|_{\mathfrak n},
	\qquad
	g\longmapsto\exp(\chi(g)D)
	\]
	have the same differential.
	Since $G$ is connected,
	they coincide.
	
	On the other hand,
	$
	\chi(K)\neq\{0\},
	$
	otherwise, $\chi$ would descend
	to a smooth function $f:G/K\too\R$ such that $df=\eta$
	This is impossible because
	$\eta$ is nowhere vanishing and $f$ has a critical point.
	
	We may therefore choose
	$k_0\in K$ such that
	$
	t_0=\chi(k_0)>0.
	$
	By \eqref{eq:period-ad},
	$
	\Ad_{k_0}|_{\mathfrak n}
	=\exp(t_0D).
	$
	Since $\Omega$ and $\theta$ are invariant under
	the isotropy action, we obtain,
	for every $Z\in\mathfrak n$,
	\[
	c(Z)
	=\Omega(Z,T)
	=\Omega(\Ad_{k_0}Z,
	\Ad_{k_0}T)
	=c(\Ad_{k_0}Z)
	\theta(\Ad_{k_0}T)
	=c(\exp(t_0D)Z).
	\]
	
	Consequently,
	$
	\exp(t_0D^*)c=c.
	$
	Since $D^*$ commutes with its
	exponential, it follows that
	\[
	\exp(t_0D^*)(D^*)^jc
	=(D^*)^jc
	\]
	for every $j\geq0$.
	By \eqref{eq:mixed-cyclic},
	these vectors span
	$\mathfrak n^*$.
	Hence
	$
	\exp(t_0D^*)=\operatorname{Id}_{\mathfrak n^*},
	$
	and therefore
	$
	\exp(t_0D)=\operatorname{Id}_{\mathfrak n}.
	$
	It follows that the subgroup
	$
	P_D=
	\{t\in\mathbb R:\exp(tD)=\operatorname{Id}_{\mathfrak n}\}
	$
	is nontrivial.
		Since $D\neq0$, $P_D$ is a
	proper closed subgroup of
	$\mathbb R$. Consequently,
	$			P_D=\tau\mathbb Z
	$
	for a uniquely determined
	$\tau>0$.
	
	\medskip
	\noindent
	\textbf{5. Dimension formula.}
	
	Recall that
	$
	\dim\mathcal O_\be
	=\dim\mathfrak s-\dim\mathfrak s_\be.
	$
	
	In Model~\textup{(A1)}, since
	$\be|_{\mathfrak s_\be}\neq0$,
	$
	\dim\kalg=\dim\mathfrak s_\be-1.
	$
	Therefore
	\[
	\dim(G/K)
	=\dim\mathfrak s+1
	-(\dim\mathfrak s_\be-1)
	=\dim\mathcal O_\be+2.
	\]
	
	In Model~\textup{(A2)}, we have
	$
	\dim\kalg=\dim\mathfrak s_\be,
	$
	and hence the same formula.
	
	Finally, in Model~\textup{(B)},
	since $c\neq0$ we can choose $Z_0\in\nil$ such that $c(Z_0)=1$. Then the linear form $\la:\mathfrak s_\be\oplus\nil\too \R$, $X+Z\mapsto\be(X)+c(Z)$ satisfies $\la(0+Z_0)=1\not=0$ and hence
	$
	\dim\kalg=\ker\la
	=\dim\mathfrak s_\be+\dim\nil-1.
	$
	Consequently,
	\[
	\begin{aligned}
		\dim(G/K)
		&=\dim\mathfrak s+1+\dim\nil
		-\dim\kalg
		=\dim\mathfrak s
		-\dim\mathfrak s_\be+2
		=\dim\mathcal O_\be+2.
	\end{aligned}
	\]
	
	The three models are mutually
	exclusive, according to whether
	$\theta$ vanishes on $\nil$
	and, when it does not, according
	to the dimension of the radical. 
	
	This completes the proof.
	\end{proof}

Compact homogeneous spaces of solvable Lie groups
may carry invariant symplectic structures
(see \cite{ZwartBoothby}). In contrast, our
structure theorem yields the following
nonexistence result for invariant LCS structures
with nonzero Lee form.

\begin{corollary}\label{cor:solvable-nonexistence}
	Let $G$ be a connected solvable Lie group
	acting transitively on a compact homogeneous
	space $M=G/K$ of dimension at least four.
	Then $M$ admits no $G$-invariant LCS
	structure with nonzero Lee form.
\end{corollary}

\begin{proof}
	Let $N$ be the kernel of the action of $G$
	on $M$. Then $G/N$ is a connected solvable
	Lie group acting effectively and transitively
	on $M$.
	
	If $M$ admitted a $G$-invariant LCS
	structure with nonzero Lee form, this
	structure would also be invariant under
	$G/N$. However,
	Theorem~\ref{thm:nf-main} asserts that
	the Lie algebra of any connected Lie group
	acting almost effectively and transitively
	on such a manifold has a nontrivial
	semisimple Levi factor.
		This contradicts the solvability of $G/N$.
\end{proof}

\begin{corollary}[Characterization of strongly
	exact invariant LCS structures]
	\label{cor:strongly-exact-reductive}
	Let $G$ be a connected Lie group acting
	almost effectively and transitively on a
	compact homogeneous space $M=G/K$ of
	dimension at least four, and let
	$(\omega,\eta)$ be a $G$-invariant LCS
	structure with nonzero Lee form.
	
	Let $(\Omega,\theta)$ be the corresponding
	infinitesimal data, and let $\mathfrak n$
	and $\mathfrak r$ denote the nilradical
	and the solvable radical of $\mathfrak g$,
	respectively.
	
	Then the following assertions are
	equivalent:
	\begin{enumerate}
		\item The LCS structure is strongly exact,
		that is, there exists a $G$-invariant
		one-form $\alpha$ on $M$ such that
		\[
		\omega=d\alpha-\eta\wedge\alpha.
		\]
		
		\item $\theta(\mathfrak n)\neq0$.
		
		\item The solvable radical of $\mathfrak g$
		is central, and
		\[
		\mathfrak g=\mathfrak s\oplus\mathbb R^q,
		\qquad
		\mathfrak r=\mathfrak n=Z(\mathfrak g),
		\qquad 1\leq q\leq2,
		\]
		where $\mathfrak s$ is a nontrivial
		compact semisimple Lie algebra.
	\end{enumerate}
\end{corollary}

\begin{proof}
	We first prove $(1)\Rightarrow(2)$.
	
	Suppose that the LCS structure is strongly
	exact. Let $\lambda\in\mathfrak g^*$ be
	the infinitesimal form corresponding
	to an invariant primitive $\alpha$.
	Then
	\[
	\Omega=d\lambda-\theta\wedge\lambda.
	\]
	By Lemma~\ref{thm:p-Boothby}, applied
	to the invariant one-form $\lambda$,
	we have
	$
	\lambda([\mathfrak n,\mathfrak g])=0.
	$
	Assume, by contradiction, that
	$\theta(\mathfrak n)=0$.
	According to Theorem~\ref{thm:nf-main},
	\[
	\mathfrak g=
	\mathfrak s\oplus
	(\mathbb RT\ltimes_D\mathfrak n),\;\theta(T)=1,\qquad
	D=\ad_T|_{\mathfrak n}\neq0\esp [\nil,\nil]=0.
	\]
	Moreover, there exists a nonzero
	$c\in\mathfrak n^*$ such that
	\[
	\Omega(Z,X)=c(Z)\theta(X),
	\qquad
	Z\in\mathfrak n,\quad X\in\mathfrak g,
	\]
	and
	$
	\ker c=\mathfrak n\cap\mathfrak k.
	$
	
	For $Z\in\mathfrak n$, strong exactness
	and the invariance of $\la$ give
	\[
	\begin{aligned}
		c(Z)
		&=\Omega(Z,T)
		=-\lambda([Z,T])
		-\theta(Z)\lambda(T)
		+\theta(T)\lambda(Z)
		=\lambda(Z).
	\end{aligned}
	\]
	Thus
	$
	c=\lambda|_{\mathfrak n}.
	$
	Consequently, for every
	$Z\in\mathfrak n$,
	\[
	c(DZ)=\lambda([T,Z])=0.
	\]
	Therefore
	$
	D(\mathfrak n)\subseteq
	\ker c\subseteq\mathfrak k.
	$
	Now $D(\mathfrak n)$ is an ideal
	of $\mathfrak g$. Indeed,
	\[
	[\mathfrak s,D(\mathfrak n)]=0,
	\;
	[\mathfrak n,D(\mathfrak n)]=0\esp [T,D(\mathfrak n)]
	=D^2(\mathfrak n)
	\subseteq D(\mathfrak n).
	\]
	Since the action is almost effective
	and $D(\mathfrak n)$ is an ideal
	contained in $\mathfrak k$, we obtain
	$
	D(\mathfrak n)=0.
	$
	Hence $D=0$, contradicting
	Theorem~\ref{thm:nf-main}.
	We conclude that
	$\theta(\mathfrak n)\neq0$.
	
	\medskip
	The implication $(2)\Rightarrow(3)$
	follows immediately from the first
	case of Theorem~\ref{thm:nf-main}.
	
	\medskip
	Finally, suppose that (3) holds.
	Then one of the two central-radical
	models of Theorem~\ref{thm:nf-main}
	occurs. We denote by $\wi\be$ the extension of $\beta$ to $\g$ by setting $\wi\be_{|\rf}=0$.
	
	In Model~\textup{(A1)}, 
	setting $\lambda=\wi\beta$, we obtain
	\[
	\Omega=d\lambda-\theta\wedge\lambda.
	\]
	
	In Model~\textup{(A2)},
	setting $\lambda=\wi\beta+z$ where $z_{|\s\oplus \R T}=0$ and $z(Z)=1$,
	we again obtain
	\[
	\Omega=d\lambda-\theta\wedge\lambda.
	\]
	In both cases, $\lambda$ vanishes
	on $\mathfrak k$ and is
	$\Ad(K)$-invariant. Indeed,
	$\beta$ is preserved by the isotropy
	action, while the radical is central.
	
	Thus $\lambda$ defines a
	$G$-invariant one-form $\alpha$
	on $M$ satisfying
	\[
	\omega=d\alpha-\eta\wedge\alpha.
	\]
	Hence the LCS structure is
	strongly exact, proving
	$(3)\Rightarrow(1)$.
\end{proof}

\begin{corollary}\label{cor:compact-lcs-lie-groups}
	Let $G$ be a compact connected Lie group of dimension $2n\geq 4$.
	Then $G$ admits a left-invariant locally conformally symplectic
	structure with nonzero Lee form if and only if $G$ is isomorphic to
	one of the following groups:
	\[
	SU(2)\times S^1,\qquad U(2),\qquad SO(3)\times S^1.
	\]
	In particular, such a group necessarily has dimension four.
\end{corollary}

\begin{proof}
	Suppose that $G$ carries a left-invariant LCS structure with
	nonzero Lee form. We apply the classification Theorem \ref{thm:nf-main} with $K=\{e\}$, so that the
	isotropy Lie algebra $\mathfrak k$ is zero. Since $G$ is compact,
	its Lie algebra is reductive with compact semisimple part:
	\[
	\mathfrak g=\mathfrak s\oplus\mathfrak z,
	\qquad \mathfrak z=Z(\mathfrak g).
	\]
	In particular, the solvable radical of $\mathfrak g$ is central.
	Therefore Model~$(B)$, which involves a nontrivial action by
	derivations, cannot occur.
	
	In Model~$(A_2)$ the isotropy Lie algebra has the form
	$
	\mathfrak k
	=\{X-\beta(X)Z:X\in\mathfrak s_\beta\},
	$
	where $\mathfrak s_\beta$ is the coadjoint stabilizer of $\beta$.
	As the displayed map from $\mathfrak s_\beta$ to $\mathfrak g$
	is injective, $\mathfrak k=0$ implies $\mathfrak s_\beta=0$.
	This is impossible for a nonzero compact semisimple Lie algebra:
	identifying $\beta$ with $b\in\mathfrak s$ by an invariant
	Euclidean inner product gives
	$\mathfrak s_\beta=Z_{\mathfrak s}(b):=\{X\in\mathfrak s,[X,b]=0\}$,
	which is nonzero since $b\in\mathfrak s_\beta$. If $\mathfrak s=0$, Model~$(A_2)$ has dimension two,
	contrary to the assumption. Thus Model~$(A_2)$ is excluded as well.
	
	Consequently, only Model~$(A_1)$ remains. In that model,
	\[
	\mathfrak g=\mathfrak s\oplus\mathbb RT,
	\qquad
	\mathfrak k=\mathfrak s_\beta\cap\ker\beta,
	\]
	where $\mathfrak s$ is compact semisimple and $\beta\ne0$.
	Write $\beta=b^\flat$ using an invariant Euclidean inner product
	on $\mathfrak s$. Then
	\[
	\mathfrak s_\beta=Z_{\mathfrak s}(b),
	\qquad
	Z_{\mathfrak s}(b)
	=\mathbb Rb\oplus
	\bigl(Z_{\mathfrak s}(b)\cap\ker\beta\bigr).
	\]
	Since $\mathfrak k=0$, it follows that
	$Z_{\mathfrak s}(b)=\mathbb Rb$. A maximal toral subalgebra
	containing $b$ lies in this centralizer, hence
	$\operatorname{rank}\mathfrak s=1$. The only nonzero compact
	semisimple Lie algebra of rank one is $\mathfrak{su}(2)$.
	Therefore
	\[
	\mathfrak g\simeq\mathfrak{su}(2)\oplus\mathbb R.
	\]
	Its simply connected Lie group is $SU(2)\times\mathbb R$.
	The classification of compact connected four-dimensional Lie groups \cite[Remark~3.23]{BiggsRemsing}, together with Corollary~\ref{cor:solvable-nonexistence}, shows that $G$ is isomorphic to one of
	$
	SU(2)\times S^1,\; U(2)$ or $SO(3)\times S^1$.

	Conversely, on $\mathfrak{su}(2)\oplus\mathbb RT$, choose
	$0\ne\beta\in\mathfrak{su}(2)^*$ and the closed one-form
	$\theta$ characterized by
	$\theta(T)=1$ and $\theta|_{\mathfrak{su}(2)}=0$.
	Model~$A_1$ yields the LCS form
	\[
	\Omega=d\beta+\beta\wedge\theta,
	\qquad d\Omega=\theta\wedge\Omega.
	\]
	It is nondegenerate, because the coadjoint stabilizer of a nonzero
	$\beta\in\mathfrak{su}(2)^*$ is one-dimensional and
	$\beta$ does not vanish on it. The corresponding left-invariant
	forms descend to every quotient of $SU(2)\times\mathbb R$ by a
	discrete central subgroup, giving the required LCS structures
	on all three compact groups.
\end{proof}

 \section{Global integration of infinitesimal models}\label{section4}
 
 We now turn to the global realization of the
 infinitesimal models obtained in
 Theorem~\ref{thm:nf-main}. Our aim is to
 determine the corresponding Lie groups and
 their isotropy subgroups.
 
 Throughout this section, let $G$ be a connected
 Lie group and $K$ a closed subgroup such that
 $M=G/K$ is compact and the action of $G$ is
 almost effective. Without loss of generality, we can suppose $G$ simply connected.  We denote by $\G$ and $\kalg$ the Lie algebras of $G$ and $K$, respectively.

 Finally, we assume that $M$ has
 dimension at least four and carries a
 $G$-invariant LCS structure $(\eta,\om)$
 with nonzero Lee form and $(\Om,\theta)\in
 \Lambda^2\G^*\times\G^*$
 is the associated infinitesimal pair.

 According to Theorem~\ref{thm:nf-main}, $\G$ is product  $\g=\s\oplus\mathfrak r$ where $\s$ is semisimple compact and $\mathfrak r$ is either central of dimension 1 or 2 or $\mathfrak r=\R T\ltimes_D\nil$. Moreover, there exists a unique $\be\in\mathfrak s^*$  determined by
 \eqref{eq:nf-beta} such that the restriction of $\be$ to any ideal of $\s$ is non zero and $\be_{|\s_\be}\not=0$ where $\s_\be=\{X\in\s,\be([X,Y])=0\;\forall\; Y\in\s\}$.

 Let $S$ be the compact connected and simply
 connected Lie group with Lie algebra
 $\mathfrak s$.
 
 The connected and simply connected Lie
 groups corresponding to the three
 infinitesimal models are
 \[
 G=
 \begin{cases}
 	S\times\mathbb R,
 	&\text{in Model (A1)},\\
 	S\times\mathbb R^2,
 	&\text{in Model (A2)},\\
 	S\times(\mathbb R\ltimes_\varphi\mathfrak n),
 	&\text{in Model (B)},
 \end{cases}
 \]
 where
 $
 \varphi:\mathbb R\longrightarrow
 GL(\mathfrak n),$ $t\mapsto
 \exp(tD).$\medskip
 
 {\bf First step: Determination of $K^0$ the connected component of $K$.}\medskip
 
 Let
 $
 H
 =\{g\in S:\Ad_g^*\beta=\beta\}.
 $
 Identifying $\beta\in\s^*$ with an element
 $b\in\mathfrak s$ by means of the
 Killing form, we have
 $H$ is the centralizer of the torus $\overline{\exp(\mathbb Rb)}$.
 Since the centralizer of a torus in
 a compact connected Lie group is
 connected (see \cite[Corollary~4.51]{Knapp}),
 $H$ is compact and connected. We denote by $\h=\s_\be$ the Lie algebra of $H$ and let
 $
 p:\widehat H\longrightarrow H
 $
 be the universal covering homomorphism.
 
 By virtue of Lemma \ref{lem:nf-compact},
 $
 \lambda=\beta|_{\mathfrak h}
 $
 is a nonzero Lie algebra character,
 since
 $
 \lambda([\mathfrak h,\mathfrak h])=0. 
 $ Consequently,  
 since  $\widehat H$ is simply connected, $\la$ integrates uniquely to a surjective Lie group homomorphism
 $
 \widehat\lambda:\widehat H\longrightarrow\mathbb R.
 $
 Moreover, $\ker\widehat\lambda$ is connected, since the kernel of any surjective Lie group homomorphism from a connected Lie group onto $(\mathbb R,+)$ is connected. Thus $\ker\widehat\lambda$ is a closed, connected Lie subgroup of $\widehat H$ with Lie algebra $\ker\lambda$.

\begin{lemma}\label{lem:global-K0} With the notations and assumptions above,
	 $K^0$ is given explicitly
	in the three models as follows.
	
	\begin{enumerate}
		\item In Model~\textup{$(A_1)$}:
		$
		K^0=
		\{(p(h),0):h\in\ker\widehat\lambda\}.
		$
		
		\item In Model~\textup{$(A_2)$}:
		$
		K^0=
		\left\{
		(p(h),0,-\widehat\lambda(h)):
		h\in\widehat H
		\right\}.
		$
		
		\item In Model~\textup{(B)}:
		$
		K^0=
		\left\{
		(p(h),0,w-\wi\la(h)v_0):
		h\in\widehat H,\;
		w\in\ker c
		\right\}
		$ where $v_0\in\nil$ is chosen such that $c(v_0)=1$.
	\end{enumerate}
\end{lemma}

\begin{proof}
	Recall that $K^0$ is the unique connected analytic
	subgroup of $G$ whose Lie algebra is $\mathfrak k$.
	According to Theorem~\ref{thm:nf-main}, we have
	\[
	\mathfrak k=
	\begin{cases}
		\mathfrak h\cap\ker\lambda,
		&\text{in Model (A1)},\\[2mm]
		\{X-\lambda(X)Z:X\in\mathfrak h\},
		&\text{in Model (A2)},\\[2mm]
		\{X+v:X\in\mathfrak h,\ v\in\mathfrak n,\
		\lambda(X)+c(v)=0\},
		&\text{in Model (B)}.
	\end{cases}
	\]
	
	In Model~\textup{(A1)}, consider the homomorphism
	\[
	\Phi_1:\ker\widehat\lambda\longrightarrow G,
	\qquad
	\Phi_1(h)=(p(h),0).
	\]
	Its image is a connected analytic subgroup with
	Lie algebra $\mathfrak h\cap\ker\lambda$.
	Hence $K^0=\Phi_1(\ker\widehat\lambda)$.
	
	In Model~\textup{(A2)}, define
	\[
	\Phi_2:\widehat H\longrightarrow G,
	\qquad
	\Phi_2(h)=(p(h),0,-\widehat\lambda(h)).
	\]
	This is a Lie group homomorphism whose image
	is connected and has Lie algebra
	\[
	\{X-\lambda(X)Z:X\in\mathfrak h\}.
	\]
	Consequently, $K^0=\Phi_2(\widehat H)$.
	
	Finally, in Model~\textup{(B)}, choose
	$v_0\in\mathfrak n$ such that $c(v_0)=1$,
	and define
	\[
	\begin{aligned}
		\Phi_3:\widehat H\times\ker c
		&\longrightarrow G,\\
		(h,w)&\longmapsto
		(p(h),0,w-\widehat\lambda(h)v_0).
	\end{aligned}
	\]
	Since the subgroup $S\times\{0\}\times\mathfrak n$
	is a direct product, $\Phi_3$ is a Lie group
	homomorphism.
	
	Its image is connected, and its Lie algebra
	consists of the elements
	\[
	X+w-\lambda(X)v_0,
	\qquad X\in\mathfrak h,\quad w\in\ker c.
	\]
	Equivalently, this Lie algebra is
	\[
	\{X+v:X\in\mathfrak h,\ v\in\mathfrak n,\
	\lambda(X)+c(v)=0\}
	=\mathfrak k.
	\]
	Thus $K^0=\Phi_3(\widehat H\times\ker c)$.
	
	This completes the proof.
\end{proof}

\medskip
\noindent
\textbf{Second step: Determination of $K$.}
\medskip

We begin with a consequence of the closedness of
$K^0$, which will play a crucial role in the
determination of $K$.

\begin{lemma}\label{lem:beta-integrality}
	Under the hypotheses of
	Theorem~\ref{thm:nf-main}, there exists a
	unique $\ell>0$ such that
	\[
	\widetilde\lambda(\ker p)=\ell\mathbb Z.
	\]
\end{lemma}

\begin{proof}
	Put
	$
	\Gamma=\ker p$ and
	$H_0=p(\ker\widetilde\lambda).
	$
	By Lemma~\ref{lem:global-K0}, we have
	$
	H_0=K^0\cap S
	$
	in all three models. Since $K^0$ and $S$
	are closed subgroups of $G$, $H_0$ is a
	closed subgroup of $H$.
	
	Moreover, $H_0$ is normal in $H$ and
	its Lie algebra is $\ker\lambda$.
	Since $\lambda\neq0$, the quotient
	$H/H_0$ is a compact connected
	one-dimensional Lie group, hence a circle.
	
	On the other hand, the surjective
	homomorphism $\widetilde\lambda$
	induces an isomorphism of topological groups
	$
	H/H_0\simeq
	\mathbb R/\widetilde\lambda(\Gamma).
	$
	Consequently, $\widetilde\lambda(\Gamma)$
	is a nontrivial discrete subgroup of
	$\mathbb R$. Thus there exists a unique
	$\ell>0$ such that
	$
	\widetilde\lambda(\ker p)=\ell\mathbb Z.
	$
\end{proof}

The preceding lemma shows that the closedness
of $K^0$ imposes a strong integrality condition
on $\lambda$. In particular, the quotient
\[
\mathbb T_\beta
:=\mathbb R/\widetilde\lambda(\ker p)
=\mathbb R/\ell\mathbb Z
\]
is a circle group.

Hence $\widetilde\lambda$ descends to a
surjective Lie group homomorphism
$
\varphi_\beta:H\longrightarrow\mathbb T_\beta,
$
uniquely determined by the commutative diagram
\[
\begin{tikzcd}[column sep=huge, row sep=huge]
	\widetilde H
	\arrow[r, "\widetilde\lambda"]
	\arrow[d, "p"']
	&
	\mathbb R
	\arrow[d, "\pi"]
	\\
	H
	\arrow[r, "\varphi_\beta"']
	&
	\mathbb T_\beta
\end{tikzcd}
\]
where $\pi$ is the canonical projection.
Identifying the Lie algebra of
$\mathbb T_\beta$ with $\mathbb R$, we have
$
(d\varphi_\beta)_e=\lambda.
$
Furthermore,
$
\ker\varphi_\beta
=p(\ker\widetilde\lambda).
$
Since $\ker\widetilde\lambda$ is connected,
so is $\ker\varphi_\beta$.

Thus $\lambda$ integrates to a surjective
character of $H$ with connected kernel.
The geometric significance of this character,
in particular its relation to integral
coadjoint orbits, will be discussed below.

We first record the adjoint actions in the three models, which will be used
in the determination of the isotropy subgroups.

\paragraph{Model~\textup{(A1)}.}
For $G_1=S\times\mathbb R$ and
$\mathfrak g_1=\mathfrak s\oplus\mathbb RT$, where $T$ is central, we have
\begin{equation}\label{eq:Ad-A1}
	\operatorname{Ad}_{(s,r)}(X+xT)=\operatorname{Ad}_sX+xT.
\end{equation}

\paragraph{Model~\textup{(A2)}.}
For $G_2=S\times\mathbb R_T\times\mathbb R_Z$ and
$\mathfrak g_2=\mathfrak s\oplus\mathbb RT\oplus\mathbb RZ$,
where $T$ and $Z$ are central, we have
\begin{equation}\label{eq:Ad-A2}
	\operatorname{Ad}_{(s,r,z)}(X+xT+yZ)
	=\operatorname{Ad}_sX+xT+yZ.
\end{equation}

\paragraph{Model~\textup{(B)}.}
Recall that $G_B=S\times(\mathbb R\ltimes_D\mathfrak n)$, where
\[
(r,u)(r',u')=(r+r',u+e^{rD}u'),
\]
and $\mathfrak g_B=\mathfrak s\oplus(\mathbb RT\ltimes_D\mathfrak n)$,
with $[T,w]=Dw$. Then
\begin{equation}\label{eq:Ad-B}
	\operatorname{Ad}_{(s,r,u)}(X+xT+w)
	=\operatorname{Ad}_sX+xT+e^{rD}w-xDu.
\end{equation}

For every integer $d\geq1$, set
\[
F_d=\left\{\left[\frac{j\ell}{d}\right]_{\ell}:
0\leq j<d\right\}\subset\mathbb T_\beta,
\qquad
H_d=\phi_\beta^{-1}(F_d)\subset H.
\]

\begin{lemma}\label{lem:global-K} Under the preceding assumptions, the closed
	isotropy subgroup $K$ is necessarily of one
	of the following forms:
	\begin{enumerate}
		\item In Model~\textup{(A1)},
		\begin{equation}\label{eq:K-explicit-A1}
			K
			=\bigl\{(h,na):h\in H,\ n\in\mathbb Z,\ 
			\phi_\beta(h)-[n\alpha]_\ell\in F_d\bigr\}.
		\end{equation}
		\item In Model~\textup{(A2)},
		\begin{equation}\label{eq:K-explicit-A2}
			K
			=\bigl\{(h,na,z):h\in H,\ n\in\mathbb Z,\ z\in\mathbb R,\ 
			\phi_\beta(h)+[z-n\alpha]_\ell\in F_d\bigr\}.
		\end{equation}
		\item In Model~\textup{(B)},
		\begin{equation}\label{eq:K-explicit-B}
			K
			=\bigl\{(h,na,u):h\in H,\ n\in\mathbb Z,\ u\in\mathfrak n,\ 
			\phi_\beta(h)+[c(u)-n\alpha]_\ell\in F_d\bigr\}.
		\end{equation}
	\end{enumerate}
	Here $d\geq1$, $a>0$, and
	$\alpha\in\mathbb R/(\ell/d)\mathbb Z$; in the displayed formulas,
	$\alpha$ denotes any real representative of this class.
	In Model~\textup{(B)}, one must also have
	$a=m\tau$ for some $m\in\mathbb Z_{>0}$, where $\tau$ is the
	fundamental period of $D$. For each fixed infinitesimal model, the
	parameters $(d,a,\alpha)$ determine $K$ uniquely.
\end{lemma}

\begin{proof}
	Let $L$ be the stabilizer of $(\Omega,\theta)$ under the adjoint
	action of $G$. Since $\theta$ is closed and $G$ is connected,
	$\theta$ is $\operatorname{Ad}(G)$-invariant. Thus $L$ is the
	stabilizer of $\Omega$; in particular, $L$ is closed in $G$ and
	contains $K$.
	
	Using Theorem~\ref{thm:nf-main} and
	\eqref{eq:Ad-A1}--\eqref{eq:Ad-B}, we obtain
	\begin{equation}\label{eq:L}
	L=
	\begin{cases}
		H\times\mathbb R,
		&\text{in Model~\textup{(A1)}},\\
		H\times\mathbb R^2,
		&\text{in Model~\textup{(A2)}},\\
		H\times\tau\mathbb Z\times\mathfrak n,
		&\text{in Model~\textup{(B)}}.
	\end{cases}
	\end{equation}
	Indeed, in Model~\textup{(B)}, invariance of $\Omega$ requires
	$c\circ e^{tD}=c$. Since $c$ is cyclic for $D^*$, this is
	equivalent to $e^{tD}=\operatorname{Id}_{\mathfrak n}$, or
	$t\in\tau\mathbb Z$. Hence the group law on $L$ in this model
	is the direct-product law.
	
	Set $Q=\mathbb T_\beta\times A$, where
	\[
	A=
	\begin{cases}
		\mathbb R,&\text{in Models~\textup{(A1)} and~\textup{(A2)}},\\
		\tau\mathbb Z,&\text{in Model~\textup{(B)}}.
	\end{cases}
	\]
	Consider the surjective Lie group homomorphism $q:L\to Q$
	defined in the three models by
	\begin{equation}\label{eq:q}
	\begin{array}{ll}
		\textup{(A1)}&q(h,t)=(\phi_\beta(h),t),\\[1mm]
		\textup{(A2)}&q(h,t,z)=(\phi_\beta(h)+[z]_\ell,t),\\[1mm]
		\textup{(B)}&q(h,t,u)=(\phi_\beta(h)+[c(u)]_\ell,t).
	\end{array}
	\end{equation}
	By Lemma~\ref{lem:global-K0}, we have $\ker q=K^0$.
	
	Put $\Delta=q(K)$. Since $q$ is a surjective Lie group
	homomorphism, it is a submersion and hence it is open. Moreover, $\ker q=K^0\subseteq K$,
	so $q^{-1}(\Delta)=K$. Consequently,
	\[
	Q\setminus\Delta=q(L\setminus K).
	\]
	The right-hand side is open, since $K$ is closed in $L$.
	Thus $\Delta$ is closed in $Q$.
	
	We next show that $\Delta$ is discrete for the topology induced
	from $Q$. Since $K^0$ is open in the Lie group $K$, there is
	an open neighborhood $U$ of the identity in $L$ such that
	$U\cap K\subseteq K^0$. As $q$ is open, $q(U)$ is an open
	neighborhood of the identity in $Q$. Furthermore,
	\[
	q(U)\cap\Delta=q(U\cap K)=\{e_Q\}.
	\]
	Hence the identity is isolated in $\Delta$, and therefore
	$\Delta$ is a closed discrete subgroup of $Q$.
	
	We now determine $\Delta$. Since $\mathbb T_\beta$ is compact,
	its projection $\operatorname{pr}_A(\Delta)$ is discrete in $A$.
	Indeed, for each compact subset $C\subseteq A$,
	\[
	\Delta\cap(\mathbb T_\beta\times C)
	\]
	is a closed discrete subset of a compact space, hence finite.
	Thus $\operatorname{pr}_A(\Delta)\cap C$ is finite.
	
	The projection $\operatorname{pr}_A(\Delta)$ is nontrivial.
	To see this, let $\chi:G\to\mathbb R$ be the surjective
	homomorphism integrating $\theta$. On $L$, the homomorphism
	$\chi$ coincides with $\operatorname{pr}_A\circ q$ (with
	$A\subseteq\mathbb R$). If $\operatorname{pr}_A(\Delta)=\{0\}$,
	then $\chi(K)=0$, so $\chi$ descends to a continuous
	surjective map $f:G/K\to\mathbb R$ such that $df=\eta$. This contradicts the
	compactness of $G/K$ and the fact that $\eta$ is nonzero everywhere. Therefore
	$
	\operatorname{pr}_A(\Delta)=a\mathbb Z
	$
	for a unique $a>0$. In Model~\textup{(B)}, necessarily
	$a=m\tau$ for some integer $m\geq1$.
	
	On the other hand, $\Delta\cap(\mathbb T_\beta\times\{0\})$
	is a closed discrete subgroup of the compact circle
	$\mathbb T_\beta$, hence is finite. Consequently,
	$
	\Delta\cap(\mathbb T_\beta\times\{0\})
	=F_d\times\{0\}
	$
	for a unique integer $d\geq1$.
	
	Choose $([\alpha]_\ell,a)\in\Delta$. For any
	$(\xi,t)\in\Delta$, there is a unique integer $n$ such that
	$t=na$. Then
	\[
	(\xi,na)-n([\alpha]_\ell,a)
	=(\xi-[n\alpha]_\ell,0)
	\in F_d\times\{0\}.
	\]
	Thus every element of $\Delta$ has a unique expression
	$(f+[n\alpha]_\ell,na)$, with $f\in F_d$ and
	$n\in\mathbb Z$, and
	\[
	\Delta=\bigl\{(f+[n\alpha]_\ell,na):
	f\in F_d,\ n\in\mathbb Z\bigr\}.
	\]
	The choice of $\alpha$ is unique modulo $\ell/d$, since any
	two elements of $\Delta$ with second component $a$ differ
	by an element of $F_d\times\{0\}$.
	
	Finally, $K=q^{-1}(\Delta)$ gives exactly the three
	expressions in the statement.
	
	 This completes the proof.
\end{proof}

Having determined the isotropy subgroups in the three
models, we now establish an important relationship
between the corresponding homogeneous spaces.
Although Models~\textup{(A1)}, \textup{(A2)}, and
\textup{(B)} involve different transitive Lie groups,
they give rise to the same LCS geometry.
More precisely, we show that the homogeneous spaces
associated with Models~\textup{(A2)} and \textup{(B)}
are LCS-isomorphic to the corresponding homogeneous
space of Model~\textup{(A1)}.
This reduction will play a key role in the proofs
of our main theorems.

\begin{proposition}[Reduction to Model~\textup{(A1)}]
	\label{prop:reduction-A1}
	Let $(G_i,K_i,\theta_i,\Omega_i)$, for $i=2,3$,
	denote the data of Models~\textup{(A2)} and
	\textup{(B)}, respectively, with parameters
	$(d,a,\alpha)$; in Model~\textup{(B)}, assume
	that $a\in\tau\mathbb Z_{>0}$.
	
	Let
	\[
	G_1=S\times\mathbb R_T,\qquad
	K_1=\bigl\{(h,na):h\in H,\ n\in\mathbb Z,\
	\phi_\beta(h)-[n\alpha]_\ell\in F_d\bigr\},
	\]
	and identify $G_1$ with the natural subgroup
	$J_i=S\times\mathbb R_T$ of $G_i$.
	
	Then
	\[
	J_i\cap K_i=K_1,\qquad G_i=J_iK_i.
	\]
	Consequently, the inclusion
	$J_i\hookrightarrow G_i$ induces a diffeomorphism
	\[
	\Phi_i:G_1/K_1\longrightarrow G_i/K_i,
	\qquad
	\Phi_i((s,t)K_1)=(s,t,0)K_i.
	\]
	
	Moreover, the invariant LCS structures induced by
	$(\theta_1,\Omega_1)$ and $(\theta_i,\Omega_i)$
	satisfy
	\[
	\Phi_i^*\eta_i=\eta_1,\qquad
	\Phi_i^*\omega_i=\omega_1.
	\]
	In particular, Models~\textup{(A2)} and
	\textup{(B)} give the same LCS homogeneous
	spaces as the corresponding instances of
	Model~\textup{(A1)}, although the transitive
	groups acting on them are different.
\end{proposition}

\begin{proof}
	The formulas for $K_2$ and $K_3$ given by \eqref{eq:K-explicit-A2} and \eqref{eq:K-explicit-B} immediately give
	$J_i\cap K_i=K_1$.
	
	Let $L_i$ be the stabilizer of
	$(\Omega_i,\theta_i)$ given by \eqref{eq:L}, and let
	$q_i:L_i\to Q_i=\mathbb T_\beta\times A_i$
	be the surjective homomorphism given by \eqref{eq:q}. Recall that
	\[
	\ker q_i=K_i^0\subseteq K_i,
	\qquad
	A_2=\mathbb R,\quad A_3=\tau\mathbb Z.
	\]
	In both models,
	\[
	J_i\cap L_i=H\times A_i,
	\qquad
	q_i(h,t,0)=(\phi_\beta(h),t).
	\]
	As $\phi_\beta:H\to\mathbb T_\beta$ is
	surjective, $q_i(J_i\cap L_i)=Q_i$.
	It follows that
	\[
	L_i=(J_i\cap L_i)\ker q_i
	\subseteq J_iK_i.
	\]
	
	On the other hand, $G_i=J_iL_i$.
	This is immediate in Model~\textup{(A2)},
	while in Model~\textup{(B)} it follows
	from the factorization
	\[
	(s,t,u)=(s,t,0)\,(e,0,e^{-tD}u),
	\]
	where $(e,0,e^{-tD}u)\in L_i$.
	Thus $G_i=J_iK_i$ and consequently,
	the natural action of $J_i$ on $G_i/K_i$
	is  transitive, with stabilizer
	$J_i\cap K_i=K_1$. Hence it induces
	the claimed diffeomorphism $\Phi_i$.
	
	Finally, the normal-form formulas give
	\[
	\theta_i\big|_{\mathfrak g_1}=\theta_1,
	\qquad
	\Omega_i\big|_{\mathfrak g_1}=\Omega_1.
	\]
	The pairs descend to invariant LCS
	structures on the corresponding quotients
	because $K_i$ stabilizes them and
	$\mathfrak k_i=\ker\Omega_i
	\subseteq\ker\theta_i$.
	
	If $\iota_i:G_1\hookrightarrow G_i$ and
	$\pi_j:G_j\to G_j/K_j$ are the natural
	maps, then
	\[
	\pi_i\circ\iota_i=\Phi_i\circ\pi_1.
	\]
	Therefore
	\[
	\begin{aligned}
		\pi_1^*(\Phi_i^*\eta_i)
		&=\iota_i^*(\pi_i^*\eta_i)
		=\theta_1^L=\pi_1^*\eta_1,\\
		\pi_1^*(\Phi_i^*\omega_i)
		&=\iota_i^*(\pi_i^*\omega_i)
		=\Omega_1^L=\pi_1^*\omega_1.
	\end{aligned}
	\]
	Since $\pi_1$ is a surjective submersion,
	its pullback on forms is injective.
	Consequently,
	\[
	\Phi_i^*\eta_i=\eta_1,\qquad
	\Phi_i^*\omega_i=\omega_1.
	\]
\end{proof}

We have now completed the global integration of the
infinitesimal models obtained in Section~\ref{section3}.
In particular, we have determined the possible
connected simply connected Lie groups $G$ and
their closed isotropy subgroups $K$ for which
$G/K$ is compact.

In the next section, we establish the converse
construction, study the geometry of the resulting
homogeneous spaces, and complete the proofs of
our main theorems.

\section{Proofs of the main theorems and geometric reduction}
\label{section5}

In this section, we prove the  main theorems
stated in the Introduction. The proofs of
Theorems~\ref{main1}, \ref{main2}, and
\ref{main3} include the geometric reduction
of the three global models obtained in
Section~\ref{section4}. Theorem~\ref{main}
then follows from these reductions and the
classification of the isotropy subgroups
established in the preceding section.

\begin{proof}[Proof of Theorem~\ref{main1}]
	Let $\beta^L$ be the left-invariant one-form on $S$ defined by
	$\beta$, and let $q_d:S\to S/H_d$ denote the quotient map.
	The subgroup $H_d=\phi_\beta^{-1}(F_d)$ is closed and normal in
	$H$, with Lie algebra $\mathfrak h\cap\ker\beta$.
	Since $H$ fixes $\beta$, the form $\beta^L$ is invariant under
	right translations by $H_d$ and vanishes on the tangent spaces
	to their orbits. Consequently, there is a unique $S$-invariant
	one-form $\alpha_\beta$ on $S/H_d$ such that
	$
	q_d^*\alpha_\beta=\beta^L.
	$
	Moreover, $\ker(d\beta^L)_e=\mathfrak h$ and
	$\beta|_{\mathfrak h}\ne0$. Hence $d\alpha_\beta$ has
	one-dimensional kernel, on which $\alpha_\beta$ is nonzero;
	thus $\alpha_\beta$ is a contact form.
	
	The formula defining $K_1$ gives a closed subgroup of
	$G_1=S\times\mathbb R$, with Lie algebra
	$
	\mathfrak k_1=\mathfrak h\cap\ker\beta.
	$
	We next verify that the pair $(\theta_1,\Omega_1)$ induces
	an invariant LCS structure on $G_1/K_1$.
	Clearly $\theta_1|_{\mathfrak k_1}=0$. Furthermore,
	\[
	\ker\Omega_1=\mathfrak h\cap\ker\beta=\mathfrak k_1.
	\]
	To see this, choose $b\in\mathfrak h$ with $\beta(b)\ne0$.
	If $Y+zT\in\ker\Omega_1$, then
	\[
	0=\Omega_1(Y+zT,b)=-z\beta(b),
	\]
	so $z=0$. The equations $\Omega_1(Y,\mathfrak s)=0$
	and $\Omega_1(Y,T)=0$ now give
	$Y\in\mathfrak h\cap\ker\beta$.
	By \eqref{eq:Ad-A1} and the $\operatorname{Ad}(H)$-invariance
	of $\beta$, both $\theta_1$ and $\Omega_1$ are
	$\operatorname{Ad}(K_1)$-invariant. Finally,
	\[
	d\theta_1=0,\qquad d\Omega_1=\theta_1\wedge\Omega_1.
	\]
	Thus they descend to a $G_1$-invariant LCS structure
	$(\eta_1,\omega_1)$ on $G_1/K_1$.
	
	Fix a real representative of the parameter $\alpha$, and choose
	$X\in\mathfrak h$ such that $\beta(X)=\alpha/a$.
	Since $(d\phi_\beta)_e=\beta|_{\mathfrak h}$, for all
	$r\in\mathbb R$ we have
	$
	\phi_\beta(\exp_H(rX))=[r\beta(X)]_\ell.
	$
	In particular, for any $n\in \mathbb Z$
	$
	(\exp_H(naX),na)\in K_1.
	$
	Consider
	\[
	\begin{aligned}
		\Psi:(S/H_d)\times(\mathbb R/a\mathbb Z)
		&\longrightarrow G_1/K_1,\\
		(sH_d,[r])&\longmapsto(s\exp_H(rX),r)K_1.
	\end{aligned}
	\]
	We first check that $\Psi$ is well defined. Suppose
	$s_1H_d=s_2H_d$ and $r_2=r_1+na$ for some
	$n\in\mathbb Z$. Put $h=s_2^{-1}s_1\in H_d$.
	Then
	\[
	\begin{aligned}
		&(s_2\exp_H(r_2X),r_2)^{-1}
		(s_1\exp_H(r_1X),r_1)
		=
		\bigl(\exp_H(-r_2X)h\exp_H(r_1X),-na\bigr).
	\end{aligned}
	\]
	The first component belongs to $H$, and its image under
	$\phi_\beta$ is
	$
	[-n\alpha]_\ell+\phi_\beta(h).
	$
	Subtracting $[-n\alpha]_\ell$ leaves
	$\phi_\beta(h)\in F_d$, so the displayed element belongs
	to $K_1$. Hence $\Psi$ is well defined.
	
	Define
	\[
	\begin{aligned}
		\Xi:G_1/K_1&\longrightarrow
		(S/H_d)\times(\mathbb R/a\mathbb Z),\\
		(s,r)K_1&\longmapsto
		\bigl(s\exp_H(-rX)H_d,[r]\bigr).
	\end{aligned}
	\]
	To check that $\Xi$ is well defined, take $(h,na)\in K_1$.
	The representatives $(s,r)$ and $(sh,r+na)$ have the same
	image under $\Xi$, because
	\[
	\phi_\beta\bigl(
	\exp_H(rX)h\exp_H(-(r+na)X)
	\bigr)
	=\phi_\beta(h)-[n\alpha]_\ell\in F_d.
	\]
	Thus the element inside $\phi_\beta$ belongs to $H_d$.
	The formulas immediately show that $\Xi$ and $\Psi$ are
	mutually inverse. Both maps are smooth, so $\Psi$ is a
	diffeomorphism. In particular, $G_1/K_1$ is compact.
	
	We now show that $\Psi$ preserves the LCS structures.
	Let $\pi_1:G_1\to G_1/K_1$ be the quotient map, and set
	\[
	\begin{aligned}
		P:S\times\mathbb R&\longrightarrow
		(S/H_d)\times(\mathbb R/a\mathbb Z),
		&(s,r)&\longmapsto(sH_d,[r]),\\
		F:S\times\mathbb R&\longrightarrow G_1,
		&(s,r)&\longmapsto(s\exp_H(rX),r).
	\end{aligned}
	\]
	By definition,
	$
	\pi_1\circ F=\Psi\circ P.
	$
	We identify forms on $S$ with their pullbacks to
	$S\times\mathbb R$. Since $(\eta_1,\omega_1)$ is induced
	by $(\theta_1,\Omega_1)$,
	\[
	\pi_1^*\eta_1=dr,\qquad
	\pi_1^*\omega_1=d\beta^L+\beta^L\wedge dr.
	\]
	For $Y\in\mathfrak s$ and $t\in\mathbb R$, a direct
	differentiation gives
	\[
	\begin{aligned}
		&T_{(s,r)}F\bigl((T_eL_s)Y,t\bigr)
		=\Bigl(
		(T_sR_{\exp_H(rX)})(T_eL_s)Y
		+t(T_eL_{s\exp_H(rX)})X,\ t\Bigr).
	\end{aligned}
	\]
	Therefore
	\[
	\begin{aligned}
		&(F^*\beta^L)_{(s,r)}\bigl((T_eL_s)Y,t\bigr)
		=\beta\bigl(\operatorname{Ad}_{\exp_H(-rX)}Y\bigr)
		+t\beta(X)
		=\beta(Y)+t\beta(X),
	\end{aligned}
	\]
	where the last equality follows from
	$\operatorname{Ad}(H)$-invariance of $\beta$. Thus
	\[
	F^*\beta^L=\beta^L+\beta(X)\,dr,
	\qquad F^*(d\beta^L)=d\beta^L.
	\]
	Using $\pi_1\circ F=\Psi\circ P$, we obtain
	\[
	\begin{aligned}
		P^*(\Psi^*\eta_1)
		&=F^*(\pi_1^*\eta_1)=dr,\\
		P^*(\Psi^*\omega_1)
		&=F^*(\pi_1^*\omega_1)\\
		&=d\beta^L+
		(\beta^L+\beta(X)\,dr)\wedge dr\\
		&=d\beta^L+\beta^L\wedge dr.
	\end{aligned}
	\]
	On the other hand, identify $\mathbb R/a\mathbb Z$
	with $S^1=\mathbb R/\mathbb Z$ via $u=[r/a]$. Since
	$q_d^*\alpha_\beta=\beta^L$, we have
	\[
	\begin{aligned}
		P^*(a\,du)&=dr,\\
		P^*(d\alpha_\beta+a\,\alpha_\beta\wedge du)
		&=d\beta^L+\beta^L\wedge dr.
	\end{aligned}
	\]
	As $P$ is a surjective submersion, pullback by $P$ is
	injective on differential forms. Consequently,
	\[
	\Psi^*\eta_1=a\,du,\qquad
	\Psi^*\omega_1=d\alpha_\beta
	+a\,\alpha_\beta\wedge du.
	\]
	Hence $\Psi$ is an LCS isomorphism.
	
	Finally, if a nonzero ideal of $\mathfrak g_1$ were contained
	in $\mathfrak k_1\subset\mathfrak s$, it would be a
	nonzero ideal of $\mathfrak s$ on which $\beta$ vanishes.
	Such an ideal contains a simple ideal, contrary to the
	assumption on $\beta$. Thus the action of $G_1$ on
	$G_1/K_1$ is almost effective.
\end{proof}

\begin{proof}[Proof of Theorem~\ref{main2}]
	
	The formula defining $K_2$ gives a closed subgroup of
	$G_2=S\times\mathbb R_T\times\mathbb R_Z$. Indeed, the
	$T$-coordinate of its elements belongs to $a\mathbb Z$,
	and for each $n\in\mathbb Z$ the condition
	\[
	\phi_\beta(h)+[z-n\alpha]_\ell\in F_d
	\]
	defines a closed subset of $H\times\mathbb R_Z$.
	Differentiating this condition at the identity gives
	\[
	\mathfrak k_2=
	\{X-\beta(X)Z:X\in\mathfrak h\}.
	\]
	
	We show that the pair
	\[
	\theta_2=t,\qquad
	\Omega_2=d_{\mathfrak s}\beta+(\beta+z)\wedge t
	\]
	induces a $G_2$-invariant LCS structure on $G_2/K_2$.
	Clearly, $\theta_2|_{\mathfrak k_2}=0$. Moreover,
	\[
	\ker\Omega_2=\mathfrak k_2.
	\]
	Indeed, let $Y+rT+vZ\in\ker\Omega_2$.
	Evaluating against an element $b\in\mathfrak h$
	such that $\beta(b)\ne0$ gives $r=0$.
	Evaluating against $\mathfrak s$ then gives
	$Y\in\mathfrak h$, while evaluation against
	$T$ yields $v=-\beta(Y)$.
	
	By \eqref{eq:Ad-A2} and the
	$\operatorname{Ad}(H)$-invariance of $\beta$,
	both $\theta_2$ and $\Omega_2$ are
	$\operatorname{Ad}(K_2)$-invariant. Since
	\[
	d\theta_2=0,\qquad
	d\Omega_2=\theta_2\wedge\Omega_2,
	\]
	they descend to a $G_2$-invariant LCS structure
	$(\eta_2,\omega_2)$ on $G_2/K_2$.
	
	The action of $G_2$ is almost effective.
	Indeed, let $I\subseteq\mathfrak k_2$ be
	an ideal of $\mathfrak g_2$, and let
	$\mathfrak s_0$ be its projection onto
	$\mathfrak s$. Since
	$
	\mathfrak g_2=
	\mathfrak s\oplus\mathbb RT\oplus\mathbb RZ,
	$
	the subspace $\mathfrak s_0$ is an ideal
	of $\mathfrak s$ contained in $\mathfrak h$.
	
	As $\mathfrak s$ is semisimple,
	$\mathfrak s_0$ is perfect. On the other hand,
	\[
	\beta([\mathfrak h,\mathfrak h])=0,
	\]
	so $\beta|_{\mathfrak s_0}=0$.
	By the assumption that $\beta$ does not vanish
	on any simple ideal of $\mathfrak s$, we obtain
	$\mathfrak s_0=0$.
	
	Since the projection
	$\mathfrak k_2\to\mathfrak h$ is injective,
	we conclude that $I=0$.
	Thus the kernel of the $G_2$-action is discrete.

	By Proposition~\ref{prop:reduction-A1}, the inclusion
	$G_1\hookrightarrow G_2$ induces an LCS isomorphism
	$
	G_1/K_1\longrightarrow G_2/K_2,
	$
	Combining this with Theorem~\ref{main1}, we conclude
	that $G_2/K_2$ is compact and LCS-isomorphic to
	$
	(S/H_d)\times S^1
	$
	endowed with the canonical LCS structure
	\[
	\eta=a\,du,\qquad
	\omega=d\alpha_\beta+a\,\alpha_\beta\wedge du.
	\]
	\end{proof}

\begin{proof}[Proof of Theorem~\ref{main3}]
	The formula defining $K_3$ gives a closed
	subgroup of $G_3$. Indeed, its $T$-coordinate
	belongs to $a\mathbb Z$, where
	$a=m\tau$ and
	$e^{\tau D}=\operatorname{Id}_{\mathfrak n}$.
	On this subgroup, the semidirect-product
	law reduces to the direct-product law.
	Moreover, for each $n\in\mathbb Z$,
	the condition
	\[
	\phi_\beta(h)+[c(u)-n\alpha]_\ell\in F_d
	\]
	defines a closed subset of $H\times\mathfrak n$.
	Differentiating at the identity gives
	\[
	\mathfrak k_3=
	\{X+u:X\in\mathfrak h,\ u\in\mathfrak n,\
	\beta(X)+c(u)=0\}.
	\]
	
	We show that the pair
	\[
	\theta_3=t,\qquad
	\Omega_3=d_{\mathfrak s}\beta+(\beta+c)\wedge t
	\]
	induces a $G_3$-invariant LCS structure
	on $G_3/K_3$.
	Clearly, $\theta_3|_{\mathfrak k_3}=0$.
	Moreover,
	\[
	\ker\Omega_3=\mathfrak k_3.
	\]
	Indeed, let $Y+rT+u\in\ker\Omega_3$.
	Evaluating against an element
	$b\in\mathfrak h$ with $\beta(b)\ne0$
	gives $r=0$.
	Evaluation against $\mathfrak s$ then gives
	$Y\in\mathfrak h$, while evaluation
	against $T$ yields
	$\beta(Y)+c(u)=0$.
	
	We next verify the
	$\operatorname{Ad}(K_3)$-invariance.
	Let $(h,na,v)\in K_3$.
	Since $na\in\tau\mathbb Z$, we have
	$e^{naD}=\operatorname{Id}_{\mathfrak n}$.
	Hence
	\[
	\operatorname{Ad}_{(h,na,v)}(Y+rT+u)
	=
	\operatorname{Ad}_hY+rT+u-rDv.
	\]
	Since $h\in H$, the form
	$d_{\mathfrak s}\beta$ is preserved.
	Furthermore, the additional terms
	in $(\beta+c)\wedge t$ arising from
	$-rDv$ cancel by skew-symmetry.
	Thus both $\theta_3$ and $\Omega_3$
	are $\operatorname{Ad}(K_3)$-invariant.
	
	Finally, since
	\[
	dc=-t\wedge D^*c,
	\]
	we have $d(c\wedge t)=0$. Consequently,
	\[
	d\theta_3=0,\qquad
	d\Omega_3=\theta_3\wedge\Omega_3.
	\]
	Therefore, the pair descends to a
	$G_3$-invariant LCS structure
	$(\eta_3,\omega_3)$ on $G_3/K_3$.
	
	It remains to show that the action is
	almost effective.
	Let $I\subseteq\mathfrak k_3$ be an ideal
	of $\mathfrak g_3$, and let
	$\mathfrak s_0$ be its projection
	onto $\mathfrak s$.
	
	As in Model~\textup{(A2)},
	$\mathfrak s_0$ is an ideal of the
	semisimple Lie algebra $\mathfrak s$
	contained in $\mathfrak h$.
	Hence it is perfect and
	$\beta|_{\mathfrak s_0}=0$.
	By our assumption on $\beta$,
	we obtain $\mathfrak s_0=0$.
	
	Therefore,
	$
	I\subseteq\mathfrak n\cap\ker c.
	$
	Since $I$ is an ideal,
	$
	D(I)=[T,I]\subseteq I.
	$
	Thus, for every $u\in I$ and every
	integer $j\geq0$, we have
	$
	c(D^ju)=0.
	$
	Since $c$ is cyclic for $D^*$,
	the covectors
	$
	c,\ D^*c,\ (D^*)^2c,\ldots
	$
	span $\mathfrak n^*$.
	Consequently, $u=0$ for every $u\in I$,
	and hence $I=0$.
	Thus the kernel of the $G_3$-action
	is discrete.

	By Proposition~\ref{prop:reduction-A1}, the inclusion
	$G_1\hookrightarrow G_i$ induces an LCS isomorphism
	$
	G_1/K_1\longrightarrow G_3/K_3,
	$
	Combining this with Theorem~\ref{main1}, we conclude
	that $G_3/K_3$ is compact and LCS-isomorphic to
	$
	(S/H_d)\times S^1
	$
	endowed with the canonical LCS structure
	\[
	\eta=a\,du,\qquad
	\omega=d\alpha_\beta+a\,\alpha_\beta\wedge du.
	\]
	\end{proof}

\begin{proof}[Proof of Theorem~\ref{main}]
		By the infinitesimal classification established
		in Section~\ref{section3}, the Lie algebra of $G$
		and its invariant LCS structure belong to one
		of Models~\textup{(A1)}, \textup{(A2)}, or
		\textup{(B)}.
		
		Since $G$ is connected and simply connected,
		these infinitesimal models determine the
		corresponding Lie groups. The classification
		of closed isotropy subgroups obtained in
		Section~\ref{section4} then shows that $K$
		has precisely one of the forms appearing in
		Theorems~\ref{main1}--\ref{main3}.
		The result follows.
	 \end{proof}
	
	\begin{proof}[Proof of Theorem \ref{thm:canonical-reduction}]
		Theorems~\ref{main1}--\ref{main3} construct the corresponding compact
	homogeneous LCS spaces and identify each of them, as an LCS manifold,
	with
	\[
	M=(S/H_d)\times S^1,\qquad
	\eta=a\,du,\qquad
	\omega=d\alpha_\beta+a\,\alpha_\beta\wedge du.
	\]
	
	The left action of $S$ on $S/H_d$ preserves $\alpha_\beta$, and
	translations of $S^1$ preserve $du$. Therefore $S\times S^1$
	acts transitively on $M$ by LCS automorphisms. The stabilizer of
	$(eH_d,[0])$ has Lie algebra
	$(\mathfrak h\cap\ker\beta)\oplus\{0\}$.
	If an ideal of $\mathfrak s\oplus\mathbb R$ is contained in this
	subalgebra, it is an ideal of $\mathfrak s$ contained in
	$\mathfrak h\cap\ker\beta$. Every nonzero ideal of the semisimple
	Lie algebra $\mathfrak s$ contains a simple ideal, contrary to the
	assumption that $\beta$ is nonzero on every simple ideal.
	Thus the action is almost effective. Moreover,
	\[
	\omega=d\alpha_\beta-\eta\wedge\alpha_\beta,
	\]
	so the LCS structure is strongly exact, with $S\times S^1$-invariant
	primitive $\alpha_\beta$.
	
	Since $H_d$ is a closed normal subgroup of $H$ and
	$H/H_d\simeq\mathbb T_\beta/F_d\simeq S^1$, the natural projection
	\[
	p:S/H_d\longrightarrow S/H
	\]
	is a principal circle bundle. We identify $S/H$ with the coadjoint
	orbit $\mathcal O_\beta$ by $sH\mapsto\operatorname{Ad}_s^*\beta$.
	Let $q:S\to S/H$ and $q_d:S\to S/H_d$ denote the quotient maps;
	thus $q=p\circ q_d$. By definition of $\alpha_\beta$,
	$
	q_d^*\alpha_\beta=\beta^L.
	$
	Recall that the Kirillov--Kostant--Souriau form \cite{Kostant1970} is defined by
	\[
	(\omega_{\mathrm{KKS}})_\mu
	(\operatorname{ad}_X^*\mu,\operatorname{ad}_Y^*\mu)
	=\mu([X,Y]),\qquad \mu\in\mathcal O_\beta.
	\]
	The form $q^*\omega_{\mathrm{KKS}}$ is left-invariant on $S$.
	At the identity, for $X,Y\in\mathfrak s$, we have
	\[
	\begin{aligned}
		(q^*\omega_{\mathrm{KKS}})_e(X,Y)
		&=(\omega_{\mathrm{KKS}})_\beta
		(\operatorname{ad}_X^*\beta,\operatorname{ad}_Y^*\beta)\\
		&=\beta([X,Y])
		=-(d\beta^L)_e(X,Y).
	\end{aligned}
	\]
	Since both forms are left-invariant, it follows that
	$
	q^*\omega_{\mathrm{KKS}}=-d\beta^L.
	$
	Consequently,
	\[
	q_d^*\bigl(p^*\omega_{\mathrm{KKS}}+d\alpha_\beta\bigr)
	=q^*\omega_{\mathrm{KKS}}+d\beta^L=0.
	\]
	As $q_d$ is a surjective submersion, its pullback on forms is
	injective. Hence
	\[
	\;d\alpha_\beta=-p^*\omega_{\mathrm{KKS}}.\;
	\]
	
	Finally, we determine the periods of $\omega_{\mathrm{KKS}}$.
	Since $S$ is simply connected and $H$ is connected, the fibration
	$H\to S\xrightarrow{q}S/H$ shows that $S/H$ is simply connected.
	Moreover, $\pi_2(S)=0$ for every Lie group. The corresponding
	homotopy exact sequence therefore yields an isomorphism
	\[
	\partial:\pi_2(S/H)\xrightarrow{\ \simeq\ }\pi_1(H).
	\]
	Let $\gamma:S^1\to H$ be a based loop. As $S$ is simply
	connected, there is a smooth map $f:D^2\to S$ with
	$f|_{\partial D^2}=\gamma$. Since $q(\gamma)$ is constant,
	$q\circ f$ induces a map $\overline f:S^2=D^2/\partial D^2\to S/H$
	representing the spherical class corresponding to $[\gamma]$.
	By Stokes' theorem,
	\[
	\begin{aligned}
		\int_{S^2}\overline f^{\,*}\omega_{\mathrm{KKS}}
		&=\int_{D^2}f^*q^*\omega_{\mathrm{KKS}}\\
		&=-\int_{D^2}f^*d\beta^L
		=-\int_\gamma\beta^L.
	\end{aligned}
	\]
	
	Let $\rho:\widetilde H\to H$ be the universal covering map,
	and let $\widetilde\lambda:\widetilde H\to\mathbb R$ be the
	character integrating $\lambda=\beta|_{\mathfrak h}$.
	If $\widetilde\gamma$ is the lift of $\gamma$ starting at the
	identity, then $\widetilde\gamma(1)\in\ker\rho$. As
	\[
	\rho^*(\beta^L|_H)=d\widetilde\lambda,
	\]
	we obtain
	\[
	\int_\gamma\beta^L
	=\widetilde\lambda\bigl(\widetilde\gamma(1)\bigr).
	\]
	Therefore the subgroup of spherical periods of
	$\omega_{\mathrm{KKS}}$ is
	\[
	-\widetilde\lambda(\ker\rho)=\ell\mathbb Z,
	\]
	where the last equality follows from
	Lemma~\ref{lem:beta-integrality}.
	
	Since $S/H$ is simply connected, the Hurewicz homomorphism
	$\pi_2(S/H)\to H_2(S/H;\mathbb Z)$ is an isomorphism.
	Thus the spherical periods are all the integral periods, and
	\[
	\left\{\int_C\omega_{\mathrm{KKS}}:
	[C]\in H_2(S/H;\mathbb Z)\right\}=\ell\mathbb Z.
	\]
	In particular, $\omega_{\mathrm{KKS}}/\ell$ has integral
	periods. By the de Rham integrality criterion,
	\[
	{\quad
		[\omega_{\mathrm{KKS}}/\ell]\in
		\operatorname{im}\!\bigl(H^2(S/H;\mathbb Z)
		\longrightarrow H^2(S/H;\mathbb R)\bigr).
		\quad}
	\]
	In fact, this integral class is primitive, since its periods
	generate $\mathbb Z$.

	Moreover, since $B=S/H$ is simply connected,
	the universal coefficient theorem implies that
	$H^2(B;\mathbb Z)$ is torsion-free. Hence the
	natural map
	\[
	H^2(B;\mathbb Z)\longrightarrow H^2(B;\mathbb R)
	\]
	is injective, proving the uniqueness of
	the integral class $u_\beta$.
		\end{proof}

The canonical description obtained above also yields
precise topological information. In particular, it
determines the low-degree homotopy and cohomology groups
and gives a uniform obstruction to the existence of
symplectic structures.

\begin{proposition}\label{thm:topological-properties}
	Let $(M,\eta,\omega)$ be a compact homogeneous LCS manifold as in
	Theorem~\ref{thm:canonical-reduction}. Write
	\[
	M\simeq P_d\times S^1,\qquad
	P_d=S/H_d,\qquad B=S/H=\mathcal O_\beta,
	\]
	where $d\geq1$, and put $r=b_2(B)$ and
	$H_0=\ker\phi_\beta$. Then:
	\begin{enumerate}
		
		\item The first homotopy groups  are
		\[
		\pi_1(M)\simeq\mathbb Z\oplus\mathbb Z/d\mathbb Z,
		\qquad \pi_2(M)\simeq\mathbb Z^{r-1}.
		\]
		
		\item The first two integral cohomology groups are
		\[
		H^1(M;\mathbb Z)\simeq\mathbb Z,
		\qquad
		H^2(M;\mathbb Z)
		\simeq\mathbb Z^{r-1}\oplus\mathbb Z/d\mathbb Z.
		\]
		In particular, $b_1(M)=1$, $b_2(M)=r-1$, and $\chi(M)=0$.
		
		\item Every class $x\in H^2(M;\mathbb R)$ satisfies
		$
		x^{\frac12\dim M}=0.
		$
		Consequently, $M$ admits no symplectic structure and, in particular,
		no K\"ahler structure.
	\end{enumerate}
\end{proposition}

\begin{proof}\begin{enumerate}
		\item We have $\pi_1(M)=\pi_1(P_d)\oplus \mathbb Z$ and $\pi_2(M)=\pi_2(P_d)$. 
		
		Recall that $H_0$ is connected and $H_d/H_0$ is isomorphic to $F_d\simeq \mathbb Z/d\mathbb Z$. 
		
		Since $\pi_1(S)=0$ and $\pi_2(S)=0$ the homotopy exact sequence of the fibration $H_d\too S\too P_d$ gives
		\[ \pi_1(P_d)\simeq\pi_0(H_d)
		\simeq\mathbb Z/d\mathbb Z\esp \pi_2(P_d)\simeq\pi_1(H_d)
		\simeq\pi_1(H_0).  \]
		To determine $\pi_1(H_0)$ let us first determine $\pi_1(H)$.

	Recall that $H$ is the centralizer of a torus
	in the compact simply connected semisimple
	Lie group $S$ and hence $H$ is connected. By the classical result on
	derived subgroups of Levi subgroups
	\cite[Proposition~12.14]{Mallet},
	the derived subgroup $H'=[H,H]$ is simply
	connected.
	Since $H/H'$ is a torus, the homotopy exact
	sequence of the fibration $H'\too H\too H/H'$ gives
	$
	\pi_1(H)\simeq\pi_1(H/H')
	\simeq\mathbb Z^m
	$
	for some integer $m\geq0$.
	
	On the other hand, since $S$ is simply
	connected, $H$ connected and $\pi_2(S)=0$, the homotopy exact sequence of the fibration
	$H\to S\to B$ gives
	\[ \pi_1(B)=0\esp 
	\pi_2(B)\simeq\pi_1(H).
	\]
	As $B$ is simply connected, the Hurewicz
	theorem and the universal coefficient theorem
	imply
	\[
	H_2(B;\mathbb Z)\simeq
	H^2(B;\mathbb Z)\simeq\mathbb Z^m.
	\]
	Consequently, $m=b_2(B)=r$.

	The fibration
	$H_0\to H\xrightarrow{\phi_\beta}\mathbb T_\beta$
	has connected fiber. Its homotopy exact sequence
	therefore gives
	\[
	0\longrightarrow\pi_1(H_0)
	\longrightarrow\pi_1(H)\longrightarrow\mathbb Z
	\longrightarrow0,
	\]
	so $\pi_1(H_0)\simeq\mathbb Z^{r-1}$.
	These identities
	prove~(1).
	\item

	Recall that, for any path-connected space $X$,
	$
	H^1(X;\mathbb Z)\simeq
	\operatorname{Hom}(\pi_1(X),\mathbb Z)
	$
	(see \cite[Section~3.1]{Hatcher}).
	Since $\pi_1(P_d)$ is finite and $B$ is
	simply connected, it follows that
	\[
	H^1(P_d;\mathbb Z)=H^1(B;\mathbb Z)=0.
	\]
	Consequently,
	the Gysin sequence of the principal circle
	bundle $p:P_d\to B$ (see \cite[Section~4.D]{Hatcher}) yields an exact sequence
	\[
	0\longrightarrow H^0(B;\mathbb Z)\simeq\mathbb Z
	\longrightarrow H^2(B;\mathbb Z)
	\xrightarrow{p^*}H^2(P_d;\mathbb Z)
	\longrightarrow0.
	\]
	As $H^2(B;\mathbb Z)\simeq\mathbb Z^r$,
	it follows that $H^2(P_d;\mathbb Z)$
	has rank $r-1$.

	By the universal coefficient theorem, the torsion
	subgroup of $H^2(P_d;\mathbb Z)$ is isomorphic
	to the torsion subgroup of $H_1(P_d;\mathbb Z)$.
	Since
	$
	H_1(P_d;\mathbb Z)\simeq\mathbb Z/d\mathbb Z,
	$
	we obtain
	$
	\operatorname{Tor}H^2(P_d;\mathbb Z)
	\simeq\mathbb Z/d\mathbb Z.
	$
	Combining this with the rank $r-1$ obtained
	from the Gysin sequence, we conclude that
	\[
	H^2(P_d;\mathbb Z)
	\simeq\mathbb Z^{r-1}\oplus\mathbb Z/d\mathbb Z.
	\]
	
	The K\"unneth formula now yields
	\[
	H^1(M;\mathbb Z)\simeq\mathbb Z,
	\qquad
	H^2(M;\mathbb Z)\simeq H^2(P_d;\mathbb Z).
	\]
	In particular, $b_1(M)=1$, $b_2(M)=r-1$,
	and $\chi(M)=0$, the last equality following
	from $M\simeq P_d\times S^1$. This proves~(2).
	
\item	
Finally, write $\dim B=2n$, so that
$\dim M=2n+2$.
Since $H^1(B;\mathbb R)=0$, the Gysin
sequence with real coefficients shows that
\[
p^*:H^2(B;\mathbb R)
\longrightarrow H^2(P_d;\mathbb R)
\]
is surjective.
Moreover, since $H^1(P_d;\mathbb R)=0$,
the K\"unneth formula shows that
\[
\operatorname{pr}_{P_d}^*:
H^2(P_d;\mathbb R)
\longrightarrow H^2(M;\mathbb R)
\]
is an isomorphism.
Consequently, every class
$x\in H^2(M;\mathbb R)$ can be written as
$
x=(p\circ\operatorname{pr}_{P_d})^*y
$
for some $y\in H^2(B;\mathbb R)$.

Since $\dim B=2n$, we have $y^{n+1}=0$.
Therefore,
$
x^{n+1}=0.
$

If $M$ admitted a symplectic form $\sigma$,
then $[\sigma]^{n+1}$ would be nonzero,
a contradiction. Hence $M$ admits no
symplectic structure and, in particular,
no K\"ahler structure. This proves (3) and completes the proof.

\end{enumerate}
\end{proof}
\begin{remark}
	The equality $b_1(M)=1$ was already obtained
	by Vaisman \cite{Vaisman} in his study of
	compact homogeneous LCS manifolds under
	additional regularity assumptions.
	
	We extend this picture by determining the
	first two homotopy groups and the low-degree
	integral cohomology. In particular, we show that
	$
	b_2(M)=b_2(\mathcal O_\beta)-1,
	$
	where $\mathcal O_\beta$ is the associated
	coadjoint orbit. Moreover, the torsion
	subgroup of $H^2(M;\mathbb Z)$ is cyclic
	of order $d$, where $d$ is one of the
	discrete parameters arising in our
	classification.
	
	We also obtain a uniform cohomological
	obstruction to the existence of symplectic
	structures: every class
	$x\in H^2(M;\mathbb R)$ satisfies
	$
	x^{\frac12\dim M}=0.
	$
	In particular, none of these manifolds
	admits a symplectic structure, even
	without any invariance assumption.
	
\end{remark}

\section{Construction of Admissible Covectors}
\label{section6}

We begin by describing a systematic procedure
for constructing the admissible covectors $\beta$
appearing in our classification.

Let $S$ be a compact, connected, simply connected
semisimple Lie group with Lie algebra $\mathfrak s$,
and fix an $\operatorname{Ad}(S)$-invariant positive
definite inner product $\langle\,,\,\rangle$.

Given $b\in\mathfrak s\setminus\{0\}$, set
\[
T_b=\overline{\{\exp_S(tb):t\in\mathbb R\}},
\qquad 
H=Z_S(T_b)
:=\{g\in S:gh=hg
\text{ for every }h\in T_b\}.
\]
Then $H$ is connected, with Lie algebra
$\mathfrak h=\ker\operatorname{ad}_b$.

Define
\[
\beta(Y)=\langle b,Y\rangle,
\qquad Y\in\mathfrak s,
\qquad
\lambda=\beta|_{\mathfrak h}.
\]
Since $b\in Z(\mathfrak h)$, we have
\[
\lambda(b)=\langle b,b\rangle>0,
\qquad
\lambda([\mathfrak h,\mathfrak h])=0.
\]

Thus the construction reduces to two questions:
determining when $\lambda$ integrates to a
 character $\phi:H\too\R/\ell\mathbb Z$ with connected kernel,
and ensuring that $\beta$ does not vanish
on any simple ideal of $\mathfrak s$.

Choose a maximal torus $T$ of $S$ containing $T_b$, let
$\mathfrak t=\operatorname{Lie}(T)$, and denote its integral lattice and its dual lattice by
\[
\Lambda_T=\ker\bigl(\exp_S:\mathfrak t\longrightarrow T\bigr)\esp \Lambda_T^\#=\{v\in\mathfrak{t}, \;\langle v,\Lambda_T\rangle\subset\mathbb Z \}.
\]
An element $v\in \Lambda_T^\#$ is called primitive if $v$ cannot be written as $k v_0$ where $k$ is an integer with $k\geq2$ and $v_0\in \Lambda_T^\#$.

Since $\lambda([\mathfrak h,\mathfrak h])=0$,
the left-invariant one-form $\lambda^L$
on $H$ is closed and defines a period
homomorphism
\[
\operatorname{Per}_{\lambda}:
\pi_1(H)\longrightarrow\mathbb R,
\qquad
[\gamma]\longmapsto
\int_\gamma\lambda^L.
\]We denote by $\mathcal{P}=\operatorname{Per}_{\lambda}(
\pi_1(H))$.

\begin{proposition}
	\label{prop:admissible-period-criterion}
	Fix $\ell>0$. The following assertions are
	equivalent:
	\begin{enumerate}
		\item There exists a surjective Lie group
		character
		$
		\phi:H\longrightarrow\mathbb R/\ell\mathbb Z
		$
		with $(d\phi)_e=\lambda$ and connected kernel.
		
		\item The period group of $\lambda$ satisfies
		$
		\mathcal P=\ell\mathbb Z.
		$
		
		\item There exists a primitive element
		$v\in\Lambda_T^\sharp$ such that
		$
		b=\ell v.
		$
	\end{enumerate}
	Moreover, $\beta$ is nonzero on every
	simple ideal of $\mathfrak s$ if and only if
	the projection of $b$ onto each simple
	ideal is nonzero.
\end{proposition}

\begin{proof} $(1)\Longrightarrow (2)$. Suppose that $(1)$ is satisfied. Let $dt$ denote the invariant one-form on
	$\mathbb R/\ell\mathbb Z$ induced by $dt$ on $\mathbb R$.
	Since both $\phi^*dt$ and $\la^L$ are left-invariant 1-forms on $H$
	and agree at the identity, they coincide. Thus, for every
	closed loop $\gamma$ in $H$,
	\[
	\int_\gamma\beta^L
	=\int_{\phi\circ\gamma}dt
	=\ell\,\deg(\phi\circ\gamma)
	\in\ell\mathbb Z.
	\]
	Consequently, $\mathcal{P}\subset\ell\mathbb Z$.

	Write $H_0=\ker\phi$. The homotopy exact sequence of the fibration
	$H_0\to H\xrightarrow{\phi}\mathbb R/\ell\mathbb Z$ ends with
	\[
	\pi_1(H)\xrightarrow{\phi_*}\pi_1(\mathbb R/\ell\mathbb Z)
	\longrightarrow\pi_0(H_0)\longrightarrow0.
	\]
	Identifying $\pi_1(\mathbb R/\ell\mathbb Z)$ with $\mathbb Z$ by winding
	number, one has
	\[
	\phi_*([\gamma])=\frac1\ell\int_\gamma\lambda^L.
	\]
	Thus $\operatorname{im}\phi_*=\mathcal P/\ell$, so exactness gives
	$
	\pi_0(H_0)\simeq\mathbb Z/(\mathcal P/\ell).
	$
	Consequently, $H_0$ is connected exactly when
	$\mathcal P=\ell\mathbb Z$.
	
	$(2)\Longrightarrow(3)$. Suppose that $(2)$ holds.

Recall that $H$ is compact and connected,
and that $T$ is a maximal torus of $H$.
By \cite[Section~0.3.3, (4)]{Wallach},
the homogeneous space $H/T$ is simply
connected. Hence the homotopy exact
sequence of
\[
T\longrightarrow H\longrightarrow H/T
\]
shows that the natural homomorphism
\[
\pi_1(T)\longrightarrow\pi_1(H)
\]
is surjective.

For every $X\in\Lambda_T$, the loop
$\gamma_X(t)=\exp_S(tX)$, $0\leq t\leq1$,
represents the corresponding element
of $\pi_1(T)$, and
\[
\int_{\gamma_X}\lambda^L
=\int_0^1\lambda(X)\,dt
=\lambda(X).
\]
By the preceding surjectivity, we conclude
that
\[
\mathcal P
=\lambda(\Lambda_T)=\{\langle b,\ga\rangle,\; \ga\in \Lambda_T\}.
\]
Hence $v=b/\ell$ belongs to
$\Lambda_T^\#$ and
$
\langle v,\Lambda_T\rangle=\mathbb Z.
$
If $v=kv_0$ for some integer $k\geq2$
and $v_0\in\Lambda_T^\#$, then
$
\langle v,\Lambda_T\rangle
\subseteq k\mathbb Z,
$
a contradiction. Thus $v$ is primitive.

$(3)\Longrightarrow(1)$. Suppose that $(3)$ holds. 
Since $v$ is primitive in the dual lattice,
its integer pairings with $\Lambda_T$
generate $\mathbb Z$. Consequently,
$
\mathcal P=\ell\mathbb Z.
$

For $h\in H$, choose a smooth path
$\gamma:[0,1]\to H$ such that
$\gamma(0)=e$ and $\gamma(1)=h$, and define
\[
\phi(h)=
\left[\int_\gamma\lambda^L\right]_\ell.
\]
This definition is independent of the
chosen path, since the difference of two
such integrals is a period belonging to
$\ell\mathbb Z$.

The closedness of $\lambda^L$ implies
that $\phi$ is smooth, while its left
invariance shows that $\phi$ is a
group homomorphism. Moreover,
$
(d\phi)_e=\lambda\neq0,
$
so $\phi$ is surjective.

Finally, since $\mathcal P=\ell\mathbb Z$,
the induced homomorphism
\[
\phi_{*}:\pi_1(H)\longrightarrow\mathbb Z
\]
is surjective. The homotopy exact sequence
of the fibration
\[
\ker\phi\longrightarrow H
\longrightarrow\mathbb R/\ell\mathbb Z
\]
then implies that $\ker\phi$ is connected.
\end{proof}

\medskip
\noindent
\textbf{Construction of admissible covectors.}
The preceding proposition provides an explicit
procedure for constructing all admissible
covectors, up to the coadjoint action of $S$.

\begin{enumerate}
	
	\item \emph{Choose the group and its integral lattice.}

	Choose a maximal torus $T\subset S$, with
	Lie algebra $\mathfrak t$, and determine
	its integral lattice
	\[
	\Lambda_T=
	\ker(\exp_S|_{\mathfrak t})
	=\mathbb Z\gamma_1\oplus\cdots
	\oplus\mathbb Z\gamma_N,
	\]
	where $N=\operatorname{rank}S$.
	
	Let $(v_1,\ldots,v_N)$ be the dual basis
	with respect to the invariant inner product:
	\[
	\langle v_i,\gamma_j\rangle=\delta_{ij}.
	\]
	Then
	\[
	\Lambda_T^\sharp=
	\mathbb Zv_1\oplus\cdots\oplus\mathbb Zv_N.
	\]
	
	\item \emph{Choose a primitive lattice vector.}
	
	Choose integers $m_1,\ldots,m_N$ such that
	\[
	\gcd(m_1,\ldots,m_N)=1,
	\]
	and put
	\[
	v=\sum_{i=1}^N m_iv_i.
	\]
	
	Write
	\[
	\mathfrak s=
	\mathfrak s_1\oplus\cdots\oplus\mathfrak s_q
	\]
	for the decomposition into simple ideals.
	Require that the orthogonal projection of
	$v$ onto every $\mathfrak s_j$ be nonzero.
	
	\item \emph{Construct the admissible covector.}
	
	Choose $\ell>0$ and define
	\[
	b=\ell v,\qquad
	\beta(Y)=\langle b,Y\rangle,
	\quad Y\in\mathfrak s.
	\]
	
	By construction, $\beta$ does not vanish
	on any simple ideal of $\mathfrak s$.
	Moreover,
	\[
	\begin{aligned}
		\be(\Lambda_T)
		&=\{\langle b,\gamma\rangle:
		\gamma\in\Lambda_T\}
		=\ell\mathbb Z.
	\end{aligned}
	\]
	
	\item \emph{Determine the stabilizer.}
	
	Compute
	\[
	H=S_{\beta}
	=Z_S(b)
	=Z_S\left(
	\overline{\{\exp_S(tb):t\in\mathbb R\}}
	\right).
	\]
	Its Lie algebra is
	$
	\mathfrak h=\ker\operatorname{ad}_b.
	$
	If $\mathcal R$ denotes the root system
	of $\mathfrak s_{\mathbb C}$ relative to
	$\mathfrak t_{\mathbb C}$, then
	\[
	\mathfrak h_{\mathbb C}
	=
	\mathfrak t_{\mathbb C}
	\oplus
	\bigoplus_{\substack{\alpha\in\mathcal R\\
			\alpha(b)=0}}
	(\mathfrak s_{\mathbb C})_\alpha.
	\]
	
	Thus $H$ can be determined explicitly
	from the roots vanishing on $b$.
	
	\item \emph{Construct the character.}
	
	Put $\lambda=\beta|_{\mathfrak h}$.
	By the preceding proposition, $\lambda$
	integrates to a surjective character
	\[
	\phi:H\longrightarrow
	\mathbb R/\ell\mathbb Z
	\]
	with connected kernel.
	Its restriction to $T$ is explicitly
	given by
	\[
	\phi(\exp_S X)
	=[\langle b,X\rangle]_\ell,
	\qquad X\in\mathfrak t.
	\]

	We have thus obtained admissible data
	$(S,\beta,\phi,\ell)$.
	For every integer $d\geq1$, define
	$
	H_d=\phi^{-1}(F_d),
	$
	as in the classification theorem.
	Choosing any $a>0$, we obtain the compact
	homogeneous LCS manifold
	$
	M=(S/H_d)\times S^1,
	$
	endowed with the canonical LCS structure
	\[
	\eta=a\,du,\qquad
	\omega=d\alpha_\beta+
	a\,\alpha_\beta\wedge du,
	\]
	where $\alpha_\beta$ is the $S$-invariant
	contact form on $S/H_d$ induced by $\beta^L$.
	
\end{enumerate}

\section{Examples}\label{section7}

We now apply the construction of admissible covectors described above.
In each example we specify the integral lattice of a maximal torus,
choose a primitive element of its dual lattice, determine the corresponding
vector $b$ and the stabilizer $H=S_{\beta_b}$, and write down the
character $\phi_{\beta_b}$ and its subgroups $H_d$.
These steps produce the homogeneous contact manifold $P_d=S/H_d$ and
hence the compact homogeneous LCS manifold
\begin{equation}\label{eq:examples-LCS}
	M_d=P_d\times S^1,\qquad
	\eta=a\,du,\qquad
	\omega=d\alpha_{\beta_b}+a\,\alpha_{\beta_b}\wedge du,
	\qquad a>0,
\end{equation}
where $S^1=\mathbb R/\mathbb Z$ and $\alpha_{\beta_b}$ is the contact
form induced by $\beta_b^L$.
We take $\ell=2\pi$ throughout; the normalization of the invariant
inner product is $\langle X,Y\rangle=-\operatorname{tr}(XY)$ on each
special unitary factor. These examples use Model~\textup{(A1)} with
$\alpha=0$; the reduction proposition also provides realizations by
Models~\textup{(A2)} and~\textup{(B)}.

\subsection{The integral lattice of the diagonal torus in $SU(N)$}

We first record a calculation that will be used repeatedly. Let
$S=SU(N)$, $N\geq2$, and let $T$ be its diagonal maximal torus. Then
\[
\mathfrak t=\{i\operatorname{diag}(\theta_1,\ldots,\theta_N):
\theta_j\in\mathbb R,\ \sum_j\theta_j=0\},
\]
and
\[
\Lambda_T=\bigoplus_{j=1}^{N-1}\mathbb Z\gamma_j,
\qquad \gamma_j=2\pi i(E_{jj}-E_{NN}).
\]
For $1\leq j\leq N-1$ put
\[
v_j=\frac{i}{2\pi}\left(E_{jj}-\frac1N I_N\right).
\]
A direct computation gives
$\langle v_j,\gamma_k\rangle=\delta_{jk}$; hence
$(v_1,\ldots,v_{N-1})$ is the dual lattice basis.
For integers $m_1,\ldots,m_{N-1}$ with
$\gcd(m_1,\ldots,m_{N-1})=1$, choose
\[
v=\sum_{j=1}^{N-1}m_jv_j,
\qquad b=2\pi v
=i\operatorname{diag}(m_1-\overline m,\ldots,
m_{N-1}-\overline m,-\overline m),
\quad \overline m=\frac1N\sum_{j=1}^{N-1}m_j.
\]
Then
\[
\langle b,\gamma_j\rangle=2\pi m_j,
\qquad \lambda_b(\Lambda_T)=2\pi\mathbb Z.
\]
In particular, the resulting character on $H=Z_{SU(N)}(b)$
has connected kernel. Its restriction to $T$ is
\begin{equation}\label{eq:examples-SUN-character}
	\phi_{\beta_b}\bigl(\operatorname{diag}(e^{i\theta_1},\ldots,
	e^{i\theta_N})\bigr)
	=\left[\sum_{j=1}^{N-1}m_j\theta_j\right]_{2\pi}.
\end{equation}
The multiplicities of the diagonal entries of $b$ determine $H$.

\subsection{Projective spaces, lens spaces, and the four-dimensional case}

\noindent\emph{Choice of the primitive vector and of $b$.}
Take $S=SU(n+1)$, $n\geq1$, and choose
\[
(m_1,\ldots,m_n)=(1,0,\ldots,0).
\]
This vector is primitive. The preceding formula gives
\[
b=\frac{i}{n+1}\operatorname{diag}(n,-1,\ldots,-1),
\qquad \beta_b(Y)=-\operatorname{tr}(bY).
\]
It is nonzero, so $\beta_b$ does not vanish on the only simple
ideal $\mathfrak{su}(n+1)$, and
$\lambda_b(\Lambda_T)=2\pi\mathbb Z$.

\noindent\emph{Stabilizer and primitive character.}
The eigenspaces of $b$ have dimensions $1$ and $n$. Hence
\[
H=Z_{SU(n+1)}(b)=S(U(1)\times U(n))
=\{\operatorname{diag}(z,A):z\in U(1),\ A\in U(n),\ z\det A=1\},
\]
with $S/H\simeq\mathbb CP^n$. By
\eqref{eq:examples-SUN-character}, the character is
\[
\phi_{\beta_b}(\operatorname{diag}(z,A))=[\arg z]_{2\pi}.
\]
Its connected kernel is
$H_0=\{\operatorname{diag}(1,A):A\in SU(n)\}\simeq SU(n)$.
For $\mu_d=\{z\in S^1:z^d=1\}$ we have
\[
H_d=\{\operatorname{diag}(z,A)\in H:z\in\mu_d\},
\qquad H_d/H_0\simeq\mu_d.
\]

\noindent\emph{The LCS manifold.}
The first-column map identifies $SU(n+1)/H_0$ with the unit
sphere $S^{2n+1}\subset\mathbb C^{n+1}$.
Under this identification, $H_d/H_0$ acts by scalar multiplication.
Thus
\[
P_d\simeq L^{2n+1}(d)=S^{2n+1}/\mu_d,
\qquad
M_d\simeq L^{2n+1}(d)\times S^1,
\]
with the LCS structure~\eqref{eq:examples-LCS}.
In particular,
\[
\pi_1(M_d)\simeq\mathbb Z\oplus\mathbb Z/d\mathbb Z,
\qquad H^2(M_d;\mathbb Z)\simeq\mathbb Z/d\mathbb Z.
\]
For $n=1$ we obtain the four-dimensional spaces
$L^3(d)\times S^1$; for $n=3$ we obtain eight-dimensional
examples with $b_2=0$.

In dimension four, the above family exhausts the possible underlying
manifolds: the associated coadjoint orbit has real dimension two,
forcing $\mathfrak s\simeq\mathfrak{su}(2)$ under our nonvanishing
hypothesis on every simple factor.

For an explicit formula when $n=1$, choose a basis
$(E_1,E_2,E_3)$ of $\mathfrak{su}(2)$ such that
$E_3=i\operatorname{diag}(1,-1)$ and
$[E_1,E_2]=2E_3$, with the analogous cyclic brackets.
If $(\sigma_1,\sigma_2,\sigma_3)$ is its dual left-invariant
coframe, then $\beta_b=\sigma_3$. On $SU(2)\times\mathbb R$
the model forms are
\[
\theta=dr,\qquad
\Omega=-2\sigma_1\wedge\sigma_2+\sigma_3\wedge dr.
\]
They descend to $L^3(d)\times S^1$ for $\alpha=0$, with
$dr=a\,du$.

Among these four-dimensional manifolds, the compact connected Lie
groups carrying left-invariant LCS structures include
$SU(2)\times S^1$, $U(2)$, and $SO(3)\times S^1$.
In Model~\textup{(A1)} they arise from the central discrete
isotropy groups
\[
\begin{array}{c|c}
	\text{isotropy in }SU(2)\times\mathbb R&\text{quotient group}\\ \hline
	\langle(I,a)\rangle&SU(2)\times S^1\\
	\langle(-I,a)\rangle&U(2)\\
	\langle(-I,0),(I,a)\rangle&SO(3)\times S^1.
\end{array}
\]
These correspond respectively to
$(d,\alpha)=(1,0)$, $(1,\pi)$, and $(2,0)$.

\subsection{A six-dimensional family from $SU(2)\times SU(2)$}

\noindent\emph{Integral lattice and primitive vector.}
Let $S=SU(2)\times SU(2)$ and let $T=T_1\times T_2$, where
$T_j=\{\operatorname{diag}(e^{i\theta_j},e^{-i\theta_j})\}$.
Its integral lattice has basis
\[
\gamma_1=(2\pi i\operatorname{diag}(1,-1),0),\qquad
\gamma_2=(0,2\pi i\operatorname{diag}(1,-1)),
\]
and the corresponding dual vectors are
\[
v_1=\left(\frac{i}{4\pi}\operatorname{diag}(1,-1),0\right),
\qquad
v_2=\left(0,\frac{i}{4\pi}\operatorname{diag}(1,-1)\right).
\]
Choose nonzero integers $p,q$ with $\gcd(p,q)=1$ and put
\[
b=2\pi(pv_1+qv_2)
=\left(\frac{ip}{2}\operatorname{diag}(1,-1),
\frac{iq}{2}\operatorname{diag}(1,-1)\right).
\]
Then
\[
\langle b,\gamma_1\rangle=2\pi p,
\qquad \langle b,\gamma_2\rangle=2\pi q,
\qquad \lambda_b(\Lambda_T)=2\pi\mathbb Z.
\]
The condition $pq\ne0$ ensures that $\beta_b$ is nonzero on
both simple factors.

\noindent\emph{Stabilizer and character.}
Since both factors of $b$ are nonzero, $H=Z_S(b)=T^2$,
and $B=S/H\simeq S^2\times S^2$.
The character and its connected kernel are
\[
\begin{aligned}
	\phi_{\beta_b}(\theta_1,\theta_2)
	&=[p\theta_1+q\theta_2]_{2\pi},\\
	H_0&=\{(\operatorname{diag}(e^{iqt},e^{-iqt}),
	\operatorname{diag}(e^{-ipt},e^{ipt})):t\in\mathbb R\}.
\end{aligned}
\]
The kernel is connected because $(p,q)$ is primitive. More generally,
\[
H_d=\{(\theta_1,\theta_2)\in T^2:
e^{id(p\theta_1+q\theta_2)}=1\}.
\]

\noindent\emph{The LCS manifold.}
Here
\[
P_d=(SU(2)\times SU(2))/H_d,
\qquad M_d=P_d\times S^1,
\]
with the structure~\eqref{eq:examples-LCS}.
Since $b_2(B)=2$, our topological proposition gives
\[
\begin{aligned}
	\pi_1(M_d)&\simeq\mathbb Z\oplus\mathbb Z/d\mathbb Z,\\
	\pi_2(M_d)&\simeq\mathbb Z,\\
	H^2(M_d;\mathbb Z)&\simeq\mathbb Z\oplus\mathbb Z/d\mathbb Z.
\end{aligned}
\]
Thus $b_2(M_d)=1$.

\subsection{Aloff--Wallach spaces from $SU(3)$}

The Aloff--Wallach spaces form a remarkable family
of compact simply connected homogeneous manifolds
of dimension seven. Introduced by Aloff and Wallach
in their study of homogeneous Riemannian manifolds
with positive sectional curvature \cite{Aloff}, they are defined by
\[
W_{k,l}=SU(3)/T^1_{k,l},
\qquad \gcd(k,l)=1,
\]
where $T^1_{k,l}$ is the one-dimensional subgroup
of diagonal matrices with weights
$(k,l,-k-l)$.

These spaces have also attracted considerable
attention because of their topology and their
special geometric structures. In particular (see \cite{Aloff}),
\[
H^4(W_{k,l};\mathbb Z)
\simeq\mathbb Z/(k^2+kl+l^2)\mathbb Z.
\]
We show that, when $kl(k+l)\neq0$,
the Aloff--Wallach spaces arise naturally
as the homogeneous contact manifolds in
our classification. Consequently, their
products with a circle carry explicit
strongly exact homogeneous LCS structures,
providing infinitely many topologically
distinct examples in dimension eight.

\noindent\emph{Integral lattice and primitive vector.}
Let $S=SU(3)$ with its diagonal maximal torus $T^2$.
Use the lattice basis
\[
\gamma_1=2\pi i\operatorname{diag}(1,0,-1),
\qquad
\gamma_2=2\pi i\operatorname{diag}(0,1,-1),
\]
and the dual vectors
\[
v_1=\frac{i}{6\pi}\operatorname{diag}(2,-1,-1),
\qquad
v_2=\frac{i}{6\pi}\operatorname{diag}(-1,2,-1).
\]
Choose coprime integers $k,l$ with $kl(k+l)\ne0$ and set
\[
v=l v_1-k v_2,
\qquad
b_{k,l}=2\pi v
=\frac{i}{3}\operatorname{diag}(k+2l,-2k-l,k-l),
\qquad
\beta_{k,l}(Y)=-\operatorname{tr}(b_{k,l}Y).
\]
Then
\[
\langle b_{k,l},\gamma_1\rangle=2\pi l,
\qquad \langle b_{k,l},\gamma_2\rangle=-2\pi k,
\]
so $v$ is primitive and $\lambda_{b_{k,l}}(\Lambda_T)=2\pi\mathbb Z$.

\noindent\emph{Stabilizer and character.}
The differences of the diagonal coefficients of $b_{k,l}$ are
$k+l$, $-k$, and $l$. Thus the hypothesis $kl(k+l)\ne0$
ensures that $b_{k,l}$ is regular and $H=T^2$.
The associated coadjoint orbit is the complete flag manifold
$B=SU(3)/T^2$ of real dimension six.
For $h=\operatorname{diag}(z_1,z_2,z_3)\in T^2$, we obtain
\[
\phi_{\beta_{k,l}}(h)=[\arg(z_1^l z_2^{-k})]_{2\pi}.
\]
Its connected kernel is the circle
\[
H_0=H_{k,l}
=\{\operatorname{diag}(e^{ikt},e^{ilt},e^{-i(k+l)t}):
t\in\mathbb R\},
\]
and
\[
H_d=\{\operatorname{diag}(z_1,z_2,z_3)\in T^2:
(z_1^l z_2^{-k})^d=1\}.
\]

\noindent\emph{The LCS manifold and its topology.}
For $d=1$, the contact manifold is the Aloff--Wallach space
\[
W_{k,l}=SU(3)/H_{k,l},
\qquad M_{k,l}=W_{k,l}\times S^1.
\]
These are eight-dimensional compact homogeneous LCS manifolds
with the structure~\eqref{eq:examples-LCS}.
Since $b_2(B)=2$, they satisfy
\[
\pi_1(M_{k,l})\simeq\mathbb Z,
\qquad b_2(M_{k,l})=1.
\]
A further invariant distinguishes infinitely many of them:
\[
H^4(W_{k,l};\mathbb Z)
\simeq\mathbb Z/(k^2+kl+l^2)\mathbb Z.
\]
Also $H^3(W_{k,l};\mathbb Z)=0$, so the integral K\"unneth
formula yields the same fourth cohomology group for
$M_{k,l}=W_{k,l}\times S^1$.
Taking $(k,l)=(m,1)$, $m\geq1$, therefore produces infinitely
many pairwise nonhomeomorphic eight-dimensional LCS manifolds.
For $d>1$ the same construction gives
$P_d=SU(3)/H_d$, a finite cyclic quotient of $W_{k,l}$.

\subsection{An eight-dimensional family with $b_2=2$}

\noindent\emph{Integral lattice and primitive vector.}
Let $S=SU(2)^3$ and choose the product of the diagonal maximal tori,
$T=T_1\times T_2\times T_3\simeq(S^1)^3$.
The integral lattice and its dual have bases
\[
\begin{aligned}
	\gamma_j&=2\pi i\operatorname{diag}(1,-1)
	\quad\text{in the $j$-th factor},\\
	v_j&=\frac{i}{4\pi}\operatorname{diag}(1,-1)
	\quad\text{in the $j$-th factor},
\end{aligned}
\qquad j=1,2,3.
\]
Choose the primitive vector $v=v_1+v_2+v_3$ and set
\[
b=2\pi v
=\left(\frac{i}{2}\operatorname{diag}(1,-1),
\frac{i}{2}\operatorname{diag}(1,-1),
\frac{i}{2}\operatorname{diag}(1,-1)\right).
\]
Then $\langle b,\gamma_j\rangle=2\pi$ for all $j$,
so the period group is $2\pi\mathbb Z$, and $b$ has a nonzero
component in every simple ideal.

\noindent\emph{Stabilizer and character.}
Every component of $b$ is nonzero, hence $H=T^3$ and
$B=S/H\simeq(S^2)^3$.
The character is
\[
\phi_{\beta_b}(\theta_1,\theta_2,\theta_3)
=[\theta_1+\theta_2+\theta_3]_{2\pi},
\]
with connected kernel
\[
H_0=\{(\theta_1,\theta_2,\theta_3)\in T^3:
e^{i(\theta_1+\theta_2+\theta_3)}=1\}\simeq T^2,
\]
and
\[
H_d=\{(\theta_1,\theta_2,\theta_3)\in T^3:
e^{id(\theta_1+\theta_2+\theta_3)}=1\}.
\]

\noindent\emph{The LCS manifold.}
The contact manifold $P_d=SU(2)^3/H_d$ has dimension seven,
so $M_d=P_d\times S^1$ has dimension eight.
Since $b_2((S^2)^3)=3$, our topological proposition gives
\[
\pi_1(M_d)\simeq\mathbb Z\oplus\mathbb Z/d\mathbb Z,
\qquad \pi_2(M_d)\simeq\mathbb Z^2,
\qquad H^2(M_d;\mathbb Z)\simeq\mathbb Z^2\oplus\mathbb Z/d\mathbb Z.
\]
Thus $b_2(M_d)=2$. Together with the seven-dimensional lens spaces
and the Aloff--Wallach examples, this realizes $b_2=0,1,2$ among
our eight-dimensional homogeneous LCS manifolds.

\subsection{Different transitive presentations}

Finally, the preceding contact--circle examples admit, in addition
to the compact transitive action of $S\times S^1$, the homogeneous
presentations of Models~\textup{(A1)} and~\textup{(A2)} and suitable
presentations of Model~\textup{(B)}.
For example, take $S=SU(2)$, $a=1$, $\alpha=0$, and any $d\geq1$.
For Model~\textup{(B)} choose
\[
V=\mathbb R^2,\quad
D=2\pi\begin{pmatrix}0&-1\\1&0\end{pmatrix},\quad
c(x,y)=x.
\]
Then $e^{D}=I$, the fundamental period is $\tau=1$,
and $c,D^*c$ span $V^*$.
The resulting three transitive groups are
\[
SU(2)\times\mathbb R,\qquad
SU(2)\times\mathbb R^2,\qquad
SU(2)\times(\mathbb R\ltimes_D\mathbb R^2),
\]
but their LCS homogeneous spaces are all LCS-isomorphic to
$L^3(d)\times S^1$ with the structure~\eqref{eq:examples-LCS}.


\begin{thebibliography}{99}
	
	\bibitem{Aloff} S. Aloff and N. R. Wallach,
	An infinite family of distinct 7-manifolds admitting positively curved
	Riemannian structures, Bull. Amer. Math. Soc. 81 (1975), 93--97.
	
	
\bibitem{angella}	D. Angella, G. Bazzonni and M. Parton,\emph{Structure of locally conformally symplectic Lie algebras
	and solvmanifolds}, Ann. Sc. Norm. Super. Pisa Cl. Sci. (5)
	Vol. XX (2020), 373-411.
	
\bibitem{bande}	G. Bande and D. Kotschick, \emph{Moser stability for locally conformally symplectic structures,}
	Proc. Amer. Math. Soc. 137 (2009), no. 7, 2419–2424
	
	\bibitem{banyaga1} A. Banyaga,\emph{ On the geometry of locally conformal symplectic manifolds}, in Infinite dimensional
	Lie groups in geometry and representation theory (Washington, DC, 2000), 79–91, World Sci.
	Publishing, River Edge, NJ.
	
	\bibitem{banyaga2} A. Banyaga, \emph{Some properties of locally conformal symplectic structures}, Comment. Math. Helv.
	77 (2002) 383–398. 
	\bibitem{banayaga3} A. Banyaga, \emph{Examples of non $d\om$-exact locally conformal symplectic structures}, J. Geom. 87
	(2007), 1–13. 
	
	
	\bibitem{BazzoniMarrero}
	G. Bazzoni and J.~C. Marrero,
	\emph{On locally conformal symplectic manifolds of the first kind},
	Bulletin des Sciences Math\'ematiques \textbf{143} (2018), 1--57.
	
	
	\bibitem{BiggsRemsing}
	R. Biggs and C. C. Remsing,
	On the classification of real four-dimensional
	Lie groups,
	\textit{J. Lie Theory} \textbf{26} (2016),
	no.~4, 1001--1035.
	
	\bibitem{Wang} W. M. Boothby and H. C. Wang, On Contact Manifolds, Annals of Mathematics, 68 (1958), no 3, 721–734.
	
	\bibitem{Borel} A. Borel, Sur la cohomologie des espaces fibres principaux et des espaces
	homogenes de groupes de Lie compacts, Ann. of Math. (2) 57 (1953), 115--207.
	
	\bibitem{Diaz} A. Díaz Miranda et A. Reventós, Homogeneous contact compact manifolds and homogeneous symplectic manifolds, Bulletin des Sciences Mathématiques, 106 (1982), 337–350.
	
	\bibitem{Duval1981}
	C. Duval,
	On the polarizers of compact semi-simple Lie groups.
	Applications,
	\textit{Ann. Inst. H. Poincar\'e Sect. A}
	\textbf{34} (1981), no.~1, 95--115.
	
	
\bibitem{guedira}	F. Guédira and A. Lichnerowicz,\emph{ Géométrie des algèbres de Lie locales de Kirillov}, J.
	Math. Pures Appl. 63 (9) (1984), 407–484.
	
\bibitem{Hatcher}	
	{Hatcher, Allen},
		\emph{Algebraic Topology},
	{Cambridge University Press},
		 {2002}.
	
	
	
	\bibitem{HilgertNeeb}
	J. Hilgert and K.-H. Neeb,
	\textit{Structure and Geometry of Lie Groups},
	Springer Monographs in Mathematics,
	Springer, New York, 2012.
	
	
	\bibitem{Knapp} {Knapp, Anthony W.},
	\emph{Lie Groups Beyond an Introduction},
	{Progress in Mathematics},
	volume     {\bf140},
	{2},
	{Birkh{\"a}user},
	{Boston},
	{2002}.
	
	
	\bibitem{Kostant1970}
	B. Kostant,
	Quantization and unitary representations,
	in \textit{Lectures in Modern Analysis and Applications III},
	Lecture Notes in Mathematics, vol.~170,
	Springer, Berlin, 1970, pp.~87--208.
	
	
\bibitem{Le}	H. V. L\'e and J. Vanzura, \emph{Cohomology theories on locally conformal symplectic manifolds,}
	Asian J. Math. 19 (2015), no. 1, 45–82.
	
\bibitem{Lee}	H. C. Lee, A kind of even-dimensional differential geometry and its application to exterior
	calculus, Amer. J. Math. 65 (1943), 433–438.
	
	

\bibitem{Mallet} {Malle, Gunter and Testerman, Donna},
	\emph{Linear Algebraic Groups and
		Finite Groups of Lie Type},
	{Cambridge Studies in Advanced Mathematics},
	volume     {\bf 133},
	 {Cambridge University Press},
	{2011}.

	
	
	
\bibitem{marrero}	J. C. Marrero, D. Martînez Torres and E. Padr\'on,\emph{ Universal models via embedding
	and reduction for locally conformal symplectic structures}, Ann. Global Anal. Geom. 40 (2011),
	no. 3, 311-337.
	
	
	
	
	\bibitem{Rawnsley}
	J. Rawnsley,
	\textit{Compact coadjoint orbits},
	arXiv:math/0306336 [math.RT], 2003.
	
	
	
	\bibitem{Vaisman}
	I. Vaisman,
	\emph{Locally conformal symplectic manifolds},
	International Journal of Mathematics and Mathematical Sciences \textbf{8} (1985), no.~3, 521--536.
	
	
	
	\bibitem{Wallach} {Wallach, Nolan R.},
		\emph{Real Reductive Groups I},
		{Pure and Applied Mathematics},
		volume    {\bf132},
		 {Academic Press},
		 {1988}.
	
	
	
	
	
	\bibitem{ZwartBoothby}
	P.~B. Zwart and W.~M. Boothby,
	\emph{On compact homogeneous symplectic manifolds},
	Annales de l'Institut Fourier \textbf{30} (1980), no.~1, 129--157.
	
	
	
	
	
	
\end{thebibliography}
\end{document}